\documentclass[a4paper,11pt]{article}
\pdfoutput=1
\usepackage{amsmath}
\usepackage[a4paper, margin=2.2cm]{geometry}

\usepackage[T1]{fontenc}

\usepackage{amsbsy}

\usepackage{graphicx}

\usepackage{amssymb,amsthm}
\usepackage{bbm}
\numberwithin{equation}{section}

\usepackage[backend=biber,url=false,style=numeric-comp,isbn=false, giveninits=true,eprint=false,sorting=none,doi=false,date=year]{biblatex}

\usepackage{hyperref}
\usepackage[nameinlink, capitalise, noabbrev]{cleveref}
\crefname{assumption}{Assumption}{Assumptions}
\crefname{equation}{}{}

\usepackage{todonotes}
\setuptodonotes{inline}
\usepackage{url}

\usepackage{subcaption}
\usepackage{xfrac}
\usepackage{nicefrac}
\usepackage{acro}
\usepackage{SAAHYKacronyms}
\usepackage[export]{adjustbox}

\makeatletter
\newcommand\appendtographicspath[1]{%
  \g@addto@macro\Ginput@path{#1}%
}
\makeatother

\newtheorem{remark}{Remark}

\newtheorem{assumption}{Assumption}
\newtheorem{definition}{Definition}
\newtheorem{theorem}{Theorem}

\numberwithin{equation}{section}

\newcommand{\N}{\mathbb{N}}

\newcommand{\R}{\mathbb{R}}

\newcommand{\norm}[1]{\left\Vert #1 \right\Vert}
\newcommand{\abs}[1]{\left\vert #1 \right\vert}

\let\span\relax
\DeclareMathOperator{\span}{span}

\newcommand{\weakto}{\rightharpoonup}

\newcommand{\calX}{\mathcal{X}}
\newcommand{\calY}{\mathcal{Y}}
\newcommand{\calZ}{\mathcal{Z}}

\newcommand{\BV}{\mathrm{BV}}
\newcommand{\TV}{\mathrm{TV}}

\newcommand{\range}[1]{\mathcal R(#1)}
\renewcommand{\ker}{\mathcal N}

\newcommand{\subdiff}{\partial}

\newcommand{\defeq}{:=}

\makeatletter
\newsavebox{\@brx}
\newcommand{\llangle}[1][]{\savebox{\@brx}{\(\m@th{#1\langle}\)}%
  \mathopen{\copy\@brx\kern-0.5\wd\@brx\usebox{\@brx}}}
\newcommand{\rrangle}[1][]{\savebox{\@brx}{\(\m@th{#1\rangle}\)}%
  \mathclose{\copy\@brx\kern-0.5\wd\@brx\usebox{\@brx}}}
\makeatother

\let\sp\relax
\newcommand{\sp}[1]{\left\langle #1 \right\rangle}

\renewcommand{\phi}{\varphi}

\newcommand{\cl}[1]{{\overline{#1}}}

\newcommand{\ind}[3]{{#1=#2,\dots,#3}}
\newcommand{\ui}{x^i}
\newcommand{\yi}{y^i}
\newcommand{\zi}{z^i}

\newcommand{\orthim}{\overline}

\newcommand{\udn}{x^\U_n}

\newcommand{\vdn}{x^\Y_n}

\newcommand{\trim}{\{x^i\}_\ind{i}{1}{n}}
\newcommand{\trdata}{\{y^i\}_\ind{i}{1}{n}}
\newcommand{\trback}{\{z^i\}_\ind{i}{1}{n}}
\newcommand{\triples}{\{x^i, y^i, z^i\}_\ind{i}{1}{n}}

\newcommand{\trimhat}{\{\orthim x^i\}_\ind{i}{1}{n}}
\newcommand{\trdatahat}{\{\orthim y^i\}_\ind{i}{1}{n}}
\newcommand{\trbackhat}{\{\orthim z^i\}_\ind{i}{1}{n}}

\newcommand{\triminf}{\{x^i\}_{i \in \N}}
\newcommand{\trdatainf}{\{y^i\}_{i \in \N}}
\newcommand{\trbackinf}{\{z^i\}_{i \in \N}}

\newcommand{\trimhatinf}{\{\orthim x^i\}_{i \in \N}}
\newcommand{\trdatahatinf}{\{\orthim y^i\}_{i \in \N}}

\newcommand{\U}{\calX}
\newcommand{\Y}{\calY}

\newcommand{\udagger}{x_{\rm true}}
\newcommand{\reg}{\mathcal J}
\newcommand{\Jminsol}{x^\dagger_\reg}
\newcommand{\qand}{\quad \text{and} \quad}
\newcommand{\unJ}{x^{n,\delta}_\reg}
\newcommand{\dJ}{\subdiff \reg}
\renewcommand{\tilde}{\widetilde}

\newcommand{\iter}[2]{{#1}^{(#2)}}
\newcommand{\trni}[2]{{#1}^{#2}}
\newcommand{\totalerror}{\xi}

\DeclareMathOperator*{\argmin}{arg\,min}
\newcommand{\TrueOp}{\ensuremath{A}}
\newcommand{\ApproxOp}{\ensuremath{\widetilde{A}}}
\newcommand{\CorrectedOp}{\ensuremath{A_\Theta}}
\newcommand{\CorrectedAd}{\ensuremath{A^*_\Phi}}
\newcommand{\ForwardCor}{\ensuremath{F_\Theta}}
\newcommand{\AdjointCor}{\ensuremath{G_\Phi}}

\newcommand{\RegularisationOp}{\ensuremath{\reg}}

\bibliography{references.bib, IP.bib}

\graphicspath{{fig_proj/},{fig_modCorr/}}

\author{Simon Arridge\thanks{Department of Computer Science, University College London, 90 High Holborn, London, WC1V 6LJ, UK. Email: \texttt{Simon.Arridge@cs.ucl.ac.uk}},  Andreas Hauptmann\thanks{Research Unit of Mathematical Sciences, University of Oulu, Pentti Kaiteran katu 1, Linnanmaa, Finland. Email: \texttt{Andreas.Hauptmann@oulu.fi}}, and  Yury Korolev\thanks{Department of Mathematical Sciences, University of Bath, Bath, BA2 7AY, UK. Email: \texttt{ymk30@bath.ac.uk}}}
\title{Inverse Problems with Learned Forward Operators}
\date{}

\begin{document}
\maketitle

\textbf{Abstract:} 
Solving inverse problems requires the knowledge of the forward operator, but accurate models can be computationally expensive and hence cheaper variants that do not compromise the reconstruction quality are desired.
This chapter reviews reconstruction methods in inverse problems with learned forward operators 
that follow two different paradigms. The first one is completely 
agnostic to the forward operator and learns its restriction to the subspace spanned by the training data. The framework of regularisation by projection is then used to find a reconstruction. The second one uses a simplified model of the physics of the measurement process and only relies on the training data to learn a model correction. We present the theory of these two approaches and compare them numerically. A common theme emerges: both methods require, or at least benefit from, training data not only for the forward operator, but also for its adjoint.



%%%%%%%%%%%%%%%%%%%%%%%%%%%%%%%%%%%%%%%%%%%%%%%%
\section{Introduction}

The quality of solutions to an inverse problem depends crucially on the availability of a reliable forward model allowing one to make accurate predictions that can be compared with measured data. Such models do not always exist due to the complexity of the phenomena involved and even when accurate models exist they may be computationally too expensive for practical use in time-critical applications. Consequently, there is a need for efficient models that allow for fast computations without sacrificing reconstruction quality. 

More efficient models can be obtained by introducing simplifying assumptions, such as neglecting scattering in X-ray imaging \cite{kak2001principles,Natterer_Wubbelling}, which can only be used in certain idealised scenarios. 
Another possibility for obtaining more efficient models is to consider coarser discretisations, for instance of the finite element mesh in PDE based models, but this may lead to a considerable loss of accuracy and hence a compensation is needed to retain sufficient reconstruction quality \cite{arridge:2006,kolehmainen2009approximation}. 
Finally, in some applications the model and solutions can be constrained to a subspace allowing for a reduced order representation of the model \cite{benner2017model,quarteroni2015reduced,dolz2021model}. 

In recent years the interest in data-driven methods has also sparked new interest in designing techniques that combine analytical and learned components in the forward model. We will start with a brief overview of some data-driven methods.

\textbf{Data projection methods. }
The idea of building a low-dimensional representation of data sets and operators between them is an established technique in statistics and forms the basis of several classical machine learning methods as well as more recent deep learning based approaches. Linear methods based on principle components analysis (PCA) and robust-PCA construct spaces such that the residual error from projection onto them is small with respect to a specified tolerance in a 2-norm or 1-norm respectively. Applying these techniques to the range and domain of an operator provides a so-called reduced order model (ROM), which can as well be applied to a PDE based operator and its inverse (the Green's operator)~\cite{dolz2021model,Quarteroni2016,Feliu-Faba2020}. Kernel-PCA provides an extension by specifying a distance between data samples through a kernel function that allows for a non-linear separation criterion between components~\cite{Vidal2005}. Independent component analysis (ICA) is another classical non-linear factorisation of data that is based on maximising the statistical independence of the estimated components~\cite{You2020}. Finally, recent developments in deep learning assume that appropriate training data lie on a manifold in an abstract latent space that is obtained by learning simultaneously an encoder to, and a decoder from, the manifold such that the  composite operator (an ``autoencoder'') minimises an appropriate loss function~\cite{Lee2020}.

Let us also mention that projections are  the basis of some classical regularisation methods~\cite{NeuSch90, plato1990regularization, Kaltenbacher:2000, Poeschl:2010, hamarik2016regularization, bredies2016least}.

%%%%%%%%
\textbf{Operator learning. } 
There exists a growing body of literature on learning mappings between infinite-dimensional spaces, referred to as ``operator learning''. One prominent example is Neural Operators~\cite{kovachki2023neural} that have a multi-layer structure similar to a conventional neural network but whose layers are infinite-dimensional operators. Examples include Fourier Neural Operators~\cite{li2020fourier} and Deep Operator Networks~\cite{lu2021learning}. 
Although neural operators are infinite-dimensional objects, they need to be discretised in practice, resulting in a conventional neural network. However, in order to ensure consistency of this discretisation with the infinite-dimensional limit, down- and upsampling operators need to be included in the architecture~\cite{bartolucci2024representation}.

While training a neural operator, just like training a finite-dimensional neural network, can be computationally expensive, random feature models (also known as kernel methods) combine optimisation with randomisation and result in a significantly simpler convex problem. Random feature models also make sense in infinite dimensions and admit efficient (Monte-Carlo) convergence rates in terms of the number of features~\cite{nelsen2021random, lanthaler2024error}. In some sense, random feature models turn the problem of learning a nonlinear operator over the input space into that of learning a linear operator over the parameter space (this is sometimes referred to as the ``kernel trick''). Learning a linear operator as an inverse problem was considered in~\cite{dehoop:2023, mollenhauer2022learning} where convergence rates have also been obtained.
Approximation rates can also be obtained for infinite-dimensional holomorphic operators, which has been done in~\cite{herrmann2022neural} based on earlier work~\cite{cohen2010convergence}.
 Barron operators are another class of infinite-dimensional operators for which efficient approximation rates have been obtained~\cite{korolev:2022}. 

 In the context of inverse problems, \cite{arndt2023invertible} proposes to learn the forward operator of an inverse problem (or its regularised inversion) based on invertible residual networks.  
Learned PDE operators have also been used in parameter estimation problems~\cite{tanyu2023deep}. 

For more details on operator learning we refer the reader to the recent review~\cite{kovachki2024operator}. 

%%%%%%%%
\textbf{Our contribution. } 
In this chapter we will review two fundamentally different paradigms for solving inverse problems with the aid of training data based either on learning a data-driven representation of the forward model following the paper \cite{asp-kor-sch-2020} or learning a correction operator to a given cheaper approximation of the forward model \cite{lunz21}. 

The first approach relies on the important observation that if the forward operator is linear then its restriction to the span of the training data can be computed \textit{without any access to the forward operator}. The method proposed in~\cite{asp-kor-sch-2020} relies on orthogonalising the training set using a Gram-Schmidt process (see also~\cite{aspri:2021data} for generalisations). While this is costly, it has to be done only once and ``offline'', i.e. before solving the actual inverse problem, and can be reused for problems with the same operator but different measurements. In some sense, this is similar to using a neural network, where training is  costly but applying a trained network is cheap. 

The second approach follows the classical model correction paradigm and assumes that a computationally inexpensive simplified model is given, such as a coarser discretisation \cite{arridge:2006} or an analytic approximation~\cite{hauptmann2018approximate}, which in itself is not sufficiently accurate to produce good reconstructions when used in a  variational reconstruction framework. Recent work in \cite{lunz21,hauptmann2018approximate} considers  learning a data-driven correction given by a neural network trained on suitable training data. We will discuss two approaches to learning such a correction either as part of a learned reconstruction operator, such as an unrolled iterative scheme \cite{hauptmann2018approximate}, or as a separate explicit correction network for the forward operator \cite{lunz21}. 

We  start with a discussion of data-driven \textit{regularisation by projection} in \cref{sec:reg-proj-top}. Most of the material is taken from~\cite{asp-kor-sch-2020}, but we also add a new discussion of the connections to iterative reconstruction methods in \cref{sec:proj-note-numerics}. Then we move on, in \cref{sec:model-correction}, to the case where an approximate model is used alongside with a learned correction. We  discuss  implicit and explicit corrections and draw connections to the previous section on regularisation by projection. Then we present  numerical experiments with these approaches in \cref{sec:PAT} for the problem of limited-view photoacoustic tomography. We briefly present  ongoing work on optimisation on learned manifolds in \cref{sec:training-trajectory} and discuss possible directions for future research in \cref{sec:conclusions}. 

We also identify a common theme: it turns out that in order to obtain good reconstructions, both methods require (or at least benefit from) training data not only for the forward operator but also for its adjoint.

%%%%%%%%%%%%%%%%%%%%%%%%%%%%%%%%%%

\subsection{Mathematical setting}
Within this chapter we will consider a linear inverse problem
\begin{equation}\label{eq:Ax=y}
    Ax=y,
\end{equation}
where $A \colon \calX \to \calY$ is a linear bounded operator acting between separable Hilbert spaces $\calX$ and $\calY$ and $A^\ast \colon \calY \to \calX$ is its adjoint. Often, the exact right-hand side in~\eqref{eq:Ax=y} is not available and we only have access to an approximation $y^\delta$ such that $\norm{y-y^\delta} \leq \delta$ for some $\delta>0$. We will describe methods for solving~\eqref{eq:Ax=y} that do not require access to the exact operator $A$ during the solution phase, but rely on training pairs/triples
\begin{equation}\label{eq:training-data}
    x^i \in \calX, \quad y^i = Ax^i\in \calY, \quad \text{and} \quad z^i = A^*Ax^i \in \calX, \quad i=1,...,n,
\end{equation}
together with (in some cases) a simplified approximate model and its adjoint  $\tilde A \colon \calX \to \calY$,  $\tilde{A}^\ast \colon \calY \to \calX$.
Borrowing terminology from tomography, we will call $x^i$'s images, $y^i$'s measurements and $z^i$'s backprojections. The collection $\triples$ from~\eqref{eq:training-data} will be referred to as \textit{training data}.

%%%%%%%%%%%%%%%%%%%%%%%%%%%%%%%%%%%%%%%%%%%%%%%%
\section{Data-driven regularisation by projection}\label{sec:reg-proj-top}

%%%%%%%%%%%%%%%%%%%%%%%%%%%%%
\subsection{Setting and main assumptions}\label{sec:reg-proj-setting}

We start with a simple but important observation that was made in~\cite{asp-kor-sch-2020}: the training data~\eqref{eq:training-data} completely describe the forward and the normal operators on the span of the training images $\span \trim$. The restriction of $A$ and $A^*A$ to this subspace can be computed using Gram-Schmidt orthogonalisation, which needs to be done only once and can be done offline, prior to solving the inverse problem. In this section we will show how such learned operators can be used for regularised inversion of~\eqref{eq:Ax=y}.

Let us first fix some notation and state our main assumptions. The exact solution of~\eqref{eq:Ax=y} (with exact forward model $A$ and noise-free measurement) will be denoted by $\udagger$. 

The spans of the training images, measurements and backprojections will be denoted by
	\begin{equation*}
	\calX_n \defeq \span\trim, \quad \calY_n \defeq \span\trdata, \quad \calZ_n \defeq \span\trback.
	\end{equation*}
Orthogonal projection operators onto $\calX_n$, $\calY_n$ and $\calZ_n$ are denoted by $P_{\calX_n}$, $P_{\calY_n}$ and $P_{\calZ_n}$, respectively.

\begin{assumption}[Independence, uniform boundedness, sequentiality, density]\label{ass_1} 
	 \ \\
	\textit{Linear independence:} For every $n \in \N$ the images $\trim$ are linearly independent.\\
	\textit{Uniform boundedness:} There exist constants $c_u,C_u>0$ such that $c_u \leq \norm{\ui} \leq C_u$  for all $i \in \N$. Hence with no loss of generality we will assume that $\norm{\ui} =1$ for all $i \in \N$. \\
	\textit{Sequentiality:} The families of training triples are nested, i.e. for every $n \in \N$ 
	\begin{equation*}%\label{ass_2} 
	\{\ui, \yi, \zi\}_\ind{i}{1}{n+1} = \triples \cup \{x^{n+1}, y^{n+1}, z^{n+1}\}.
	\end{equation*}
    Consequently, the subspaces $\calX_n$, $\calY_n$ and $\calZ_n$ are nested, that is 
	\begin{equation*}
	\calX_n \subset \calX_{n+1}, \quad \calY_n \subset \calY_{n+1}, \quad \calZ_n \subset \calZ_{n+1} \quad \text{for all $n$.}
	\end{equation*}
    \textit{Density:} the subspaces spanned by the images $\triminf$ are dense in $\calX$, that is 
	\begin{equation*}
	\cl{\bigcup_{n \in \N} \calX_n} = \calX. 
	\end{equation*}
    Consequently, the subspaces spanned by the training measurements $\trdatainf$ and backprojections $\trbackinf$ are dense in the closures of the ranges $\cl{\range{A}}$ and $\cl{\range{A^*}}$, respectively. (The last statement follows from the fact that $\cl{\range{A^*A}}=\cl{\range{A^*}}$, which is easily checked.)
\end{assumption}

%%%%%%
\subsection{The injective case: regularisation by projection}\label{sec:reg-proj}

Let $y \in \range{A}$ be the exact, noise-free right-hand side in~\eqref{eq:Ax=y} and consider the following projected problem
\begin{equation} \label{AP_n u = y}
AP_{\U_n} x = y.
\end{equation}
Its minimum-norm solution is given by
\begin{equation}\label{eq:udn1}
\udn = (AP_{\U_n})^\dagger y,
\end{equation}
where $(AP_{\U_n})^\dagger$ denotes the Moore-Penrose inverse of $AP_{\U_n}$. The superscript $^{\U}$ in $\udn$ reflects the fact that the projection in~\eqref{AP_n u = y} takes place in $\U$. 

The following result shows that in the injective case, a simple formula for $(AP_{\U_n})^\dagger$ exists.
\begin{theorem}[{\cite[Thm. 4]{asp-kor-sch-2020}}]\label{thm:pseudoinv_image}
	Let  $A$ be injective. Then the Moore-Penrose inverse of $AP_{\U_n}$ is given by
	\begin{equation*}
	(AP_{\U_n})^\dagger = A^{-1} P_{\Y_n}.
	\end{equation*}
\end{theorem}

Combining this with~\eqref{eq:udn1}, we get a simple reconstruction formula 
\begin{equation}\label{eq:mp}
	\udn = A^{-1}P_{\Y_n} y.
\end{equation}
Since $A$ is injective, the restriction of its inverse to $\calY_n$ can be computed using only the training data~\eqref{eq:training-data}.  Gram-Schmidt orthogonalisation can be used for this purpose; we refer to~\cite[Sect. 3.1]{asp-kor-sch-2020} for details.

The approach~\eqref{AP_n u = y} is known as \textit{regularisation by projection}~\cite{engl:1996}. In the model-based setting, projections are taken onto subspaces spanned by a certain number of basis functions, for example, finite elements. If the basis of singular vectors of the forward operator $A$ is used, the method reduces to the truncated singular value decomposition~\cite{engl:1996}. In our case, projections are taken onto subspaces given by training data.

%%%%%%%%%%%
\subsubsection{Convergence analysis}\label{sec:proj-conv-an}

Examples of non-convergence exist~\cite{Seidman:1980} that show that minimum-norm solutions~\eqref{eq:udn1} can diverge as $n \to \infty$ even for a noiseless measurement $y$. Therefore, without additional assumptions on the subspaces $\calX_n$, we cannot expect convergence of the reconstructions~\eqref{eq:mp}. In~\cite{asp-kor-sch-2020} sufficient conditions have been obtained that rely on the interplay between the training images $\triminf$, the exact solution $\udagger$ and the forward operator $A$.

To state these conditions, let us apply Gram-Schmidt orthogonalisation to the sequence $\triminf$, 
  obtaining an orthonormal basis $\trimhatinf$ of $\U$. Transforming accordingly the training measurements $\trdatainf$, we obtain a corresponding sequence $\trdatahatinf$ such that $A \orthim x^i = \orthim y^i$. 
 Expanding the exact solution in the basis $\trimhatinf$, we get 
\begin{equation}\label{eq:exp_u_dagger}
\udagger = \sum_{i=1}^\infty \sp{\udagger, \orthim x^i} \orthim x^i.
\end{equation}

We are now ready to state our assumptions.
\begin{assumption}\label{ass:l1_coefs_gt}
Coefficients of the expansion~\eqref{eq:exp_u_dagger} are in $\ell^1$, i.e.
\begin{equation*}
\sum_{i=1}^\infty \abs{\sp{\udagger, \orthim x^i}} < \infty.
\end{equation*}
\end{assumption}

Summability conditions such as \cref{ass:l1_coefs_gt} are common in machine learning and are used to define so-called~\textit{variation norm spaces}~\cite{barron:2008, bach:2017}. 

\begin{remark}
In~\cite{arndt2023invertible}, the authors introduce the so-called \textit{local approximation property} (Theorem 3.1) which requires that the neural network achieves certain approximation rates in the vicinity of the ground truth, but not globally. This turns out to be sufficient for the convergence of regularised solutions. 
On the conceptual level, this assumption is similar to our \cref{ass:l1_coefs_gt} which only requires a certain approximation rate by the subspaces spanned by the training images for the ground truth, but not globally on the whole space.
\end{remark}

\begin{assumption}\label{ass:l2_coefs_proj}
For every $n \in \N$ and any $i \geq n+1$ consider the following expansion of $P_{\Y_n} \orthim y^i \in \Y_n$
\begin{equation}\label{eq:qn}
P_{\Y_n} \orthim y^i = \sum_{j=1}^n \beta_j^{i,n} \orthim y^j.
\end{equation}
We assume that for every $n \in \N$ 
\begin{equation}\label{eq:assump4}
\sum_{j=1}^n (\beta_j^{i,n})^2 \leq C,\qquad \textrm{for every}\,\, i \geq n+1,
\end{equation}
where $C>0$ is a constant independent of $i$ and $n$. 
\end{assumption}
We emphasise that the expansion coefficients $\beta_j^{i,n}$ change with $n$ because $\trdatahat$ is not an orthogonal basis. 

\cref{ass:l2_coefs_proj} is far less interpretable than \cref{ass:l1_coefs_gt}. It depends on the interplay between the training images $\triminf$ and the forward operator $A$. As discussed in~\cite[Sect. 6.1.3]{asp-kor-sch-2020}, checking this assumption numerically is also problematic because computing the coefficients $\beta_j^{i,n}$ would involve inverting an ill-conditioned matrix. (Note, however, that this inversion is not required for finding the solution~\eqref{eq:mp}.) In the non-convergence example from~\cite{Seidman:1980}, \cref{ass:l2_coefs_proj} can be checked analytically (and is valid).

The above two assumptions allow us to prove uniform boundedness of minimum-norm solutions~\eqref{eq:mp}.
\begin{theorem}[{\cite[Thm. 11]{asp-kor-sch-2020}}]\label{thm:conv_lsq}
Let Assumptions~\ref{ass:l1_coefs_gt} and~\ref{ass:l2_coefs_proj} be satisfied. Then $\udn$ as defined in~\eqref{eq:mp} is uniformly bounded with respect to $n$. Consequently, we have that $\udn \weakto \udagger$ weakly along a subsequence.
\end{theorem}

Under additional assumptions it is possible to prove strong convergence~\cite[Thm. 15]{asp-kor-sch-2020}.

So far, we have considered the case of a noise-free measurement $y$ in~\eqref{eq:Ax=y}. If the measurement is noisy ($y_\delta$ such that $\norm{y-y_\delta} \leq \delta$ for some known noise level $\delta > 0$) then the parameter $n$ in the projected equation~\eqref{AP_n u = y} becomes a regularisation parameter~\cite{engl:1996} that needs to be chosen as a function of the measurement noise, the larger the noise level $\delta$ the smaller $n$. Details about our specific setting can be found in~\cite[Thm. 17]{asp-kor-sch-2020}. 

It may seem counter-intuitive that increasing the size of the training set should lead to instabilities, but in our case the parameter $n$ also controls model complexity, i.e. the number of components in the solution. By the nature of the reconstruction formula~\eqref{eq:mp}, we are in the regime where the number of parameters matches the number of data (i.e., we are not in the overparametrised regime) and hence the complexity of a model has to be controlled by the noise in the data. This is in line with classical results on training neural networks from noisy data~\cite{Burger_Engl:2000}.

%%%%%%%%%%%
\subsubsection{Dual least squares}\label{sec:dual-lsq}

Although projecting the equation~\eqref{eq:Ax=y} in  the space $\U$  as in~\eqref{AP_n u = y} does not yield a convergent solution in general, it is known that projecting~\eqref{eq:Ax=y} in the space $\Y$ yields convergent solutions. This method is also referred to as \textit{dual least squares}~\cite{engl:1996}. 

The dual least squares method consists in finding the minimum norm solution of the following problem
\begin{equation}\label{eq:Q_nAu=Q_ny}
P_{\Y_n} A x = P_{\Y_n} y,
\end{equation}
We denote the minimum norm solution of~\eqref{eq:Q_nAu=Q_ny} by $\vdn$, where the superscript $^\Y$ emphasises the fact that the projection in~\eqref{eq:Q_nAu=Q_ny} takes place in $\Y$.

The following classical result shows that $\vdn$ converges strongly to the exact solution $\udagger$ as $n \to \infty$.
\begin{theorem}[{\cite[Thm. 3.24]{engl:1996}}] \label{thm:engl_dual}
Let $y$ be the exact data in~\eqref{eq:Ax=y}. Then the minimum norm solution of~\eqref{eq:Q_nAu=Q_ny} is given by 
\begin{equation}\label{eq:vdn1}
\vdn = P_{A^*\Y_n} \udagger,
\end{equation} 
where $P_{A^*\Y_n}$ is the orthogonal projector onto the subspace $A^* \Y_n$. Consequently,
\begin{equation*}
\vdn \to \udagger \quad \text{as $n \to \infty$}.
\end{equation*}
\end{theorem}

The following result gives a simple characterisation of the Moore-Penrose inverse of $P_{\Y_n}A$, similarly to Theorem~\ref{thm:pseudoinv_image}.
\begin{theorem}[{\cite[Thm. 19]{asp-kor-sch-2020}}]\label{thm:pseudoinv_data}
Let $A$ have a dense range. Then the Moore-Penrose inverse of $P_{\Y_n}A$ is given by
\begin{equation*}
(P_{\Y_n}A)^\dagger = P_{A^*\Y_n} A^{-1}.
\end{equation*}
Hence, the minimum norm solution $\vdn$ of~\eqref{eq:Q_nAu=Q_ny} is given by
\begin{equation}\label{eq:vdn2}
\vdn = P_{A^*\Y_n} A^{-1} P_{\Y_n} y = P_{A^*\Y_n} \udn,
\end{equation}
where $\udn$ is the minimum norm solution of~\eqref{AP_n u = y} as defined in~\eqref{eq:mp}.
\end{theorem}

The space $A^*\calY_n$ is, in fact, nothing but the span of the training backprojections $\trback$
\begin{equation*}
    A^*\calY_n = \calZ_n.
\end{equation*}
Therefore, in order to compute the stable reconstruction~\eqref{eq:vdn2}, having training data for the forward operator $A$ is not sufficient, one also needs training data for the adjoint $A^*$. The need for training data for the adjoint is a topic that we will also encounter later in \cref{sec:model-correction} when we will discuss data-driven model corrections.

As in \cref{sec:proj-conv-an}, if the measurement in~\eqref{eq:Ax=y} is noisy, the model complexity (that is, the dimension of the space $n$) has to be controlled by the amount of noise in the measurement. Convergence analysis of the dual least squares method can be found in~\cite[Thm. 3.26]{engl:1996}.

For numerical experiments with data-driven regularisation by projection~\eqref{AP_n u = y} and dual least squares~\eqref{eq:Q_nAu=Q_ny} we refer to~\cite[Sec. 6]{asp-kor-sch-2020}.

%%%%%%
\subsection{The non-injective case: variational regularisation}\label{sec:proj-var-reg}

If the forward operator $A$ is not injective, the results of \cref{sec:reg-proj} do not apply because~\eqref{eq:mp} requires us to be able to apply the inverse $A^{-1}$ to elements in the span of the training measurements $\calY_n=\span\trdata$. However, the training data~\eqref{eq:training-data} still allow us to evaluate the forward operator on the span of the training images $\calX_n=\span\trim$. Indeed, using the orthonormalised system $\trimhat$ and the corresponding transformed measurements $\trdatahat$, we get that for any $x \in \calX_n$
\begin{equation*}
    x = \sum_{i=1}^n \sp{x, \orthim x^i} \orthim x^i \quad \text{and} \quad Ax = \sum_{i=1}^n \sp{x, \orthim x^i} A \orthim x^i = \sum_{i=1}^n \sp{x, \orthim x^i} \orthim y^i, \quad x \in \calX_n.
\end{equation*}
For an arbitrary $x \in \calX$, therefore, we can evaluate the restriction $AP_{\U_n}$ without having numerical access to $A$:
\begin{equation}\label{eq:APnx}
    AP_{\U_n}x = \sum_{i=1}^n  \sp{x, \orthim x^i} \orthim y^i, \quad x \in \calX.
\end{equation}

The operators $AP_{\U_n}$ approximate $A$ pointwise as $n \to \infty$; if $A$ is compact then approximation also holds in the operator norm~\cite{conway:1985}. Hence, we are in the framework of inverse problems with operator errors (e.g.,~\cite{NeuSch90,TGSYag, Poeschl:2010}). 

In this section we will study the following variational regularisation problem
\begin{equation}\label{eq:proj_var_noisy}
\min_{x \in \U} \frac12 \norm{AP_{\U_n} x - y^\delta   }^2 + \alpha \reg(x),
\end{equation}
where $\reg \colon \U \to \R_+ \cup \{+\infty\}$ is a regulariser and $\alpha>0$ a regularisation parameter. 
Before we proceed with the analysis, let us say how our setting differs from the literature.

Firstly, although the forward operator in~\eqref{eq:proj_var_noisy} is evaluated on a finite-dimensional subspace, the solution will not be finite-dimensional, in general. This is different from the setting of discretised variational regularisation~\cite{NeuSch90, Poeschl:2010}, where the solution is constrained to lie in an a priori prescribed finite-dimensional space. This is also in contrast with \cref{sec:reg-proj}, where the reconstructions~\eqref{eq:mp} and~\eqref{eq:vdn2} are linear combinations of a finite number of training points.

Secondly, classical theory of regularisation under operator errors deals with bounds in the operator norm $h_n$ such that $\norm{A-AP_{\U_n}} \leq h_n$. This is a global estimate that depends on how well the subspaces $\U_n$ agree with the operator $A$ (the ideal choice would be, obviously, the eigenspaces of $A$ corresponding to  $n$ largest eigenvalues). From the data-driven point of view, we would like to work with a local error estimate such as $\norm{(A-AP_{\U_n})\Jminsol}$ or $\norm{(I-P_{\U_n})\Jminsol}$, where  $\Jminsol$ is the $\reg$-minimising solution of~\eqref{eq:Ax=y}. Even if the global approximation error $\norm{A-AP_{\U_n}}$ is large, convergence can still be fast if the training data~\eqref{eq:training-data} are chosen well for a particular solution $\Jminsol$. This will be formalised in \cref{thm:conv_rates_var} below. Such local approximation conditions appear in other contexts as well, such as regularisation by invertible residual networks~\cite{arndt2023invertible}, as discussed earlier.

%%%%%%%%%%%%%%%%%%%%%%%%%%%
\subsubsection{Convergence analysis}

We will make the following standard assumptions. 
\begin{assumption}\label{ass:J}
	The regularisation functional $\reg \colon \U \to \R_+ \cup \{+\infty\}$ is proper, convex, lower-semicontinuous, and absolutely $p$-homogeneous ($p \geq 1$). 
\end{assumption}
Denote by $\ker(\reg)$ the kernel (zero-level set) of the regulariser $\reg$, which is a linear subspace because $\reg$ is convex and absolutely $p$-homogeneous.
\begin{assumption}\label{ass:P_n_ker_J}
 The kernel $\ker(\reg)$ satisfies $\dim(\ker(\reg)) < +\infty$ and $\reg$ is coercive on the quotient space $\calX / \ker(\reg)$. Furthermore, for all $n \in \N$ we have
 \begin{equation*}
     \ker(AP_{\U_n}) \bigcap \ker(\reg) = \{0\}.
 \end{equation*}
\end{assumption}
If $\reg$ is the Total Variation regulariser~\cite{ROF}, \cref{ass:J,ass:P_n_ker_J} are satisfied if $AP_{\U_n}:L^2\to L^2$ does not annihilate constant functions.

Existence of minimisers in~\eqref{eq:proj_var_noisy} follows from standard arguments. Convergence as $\delta \to 0$ can be ensured under the usual parameter choice rule $\alpha=\alpha(\delta,n)$.
\begin{theorem}[{\cite[Thm. 23]{asp-kor-sch-2020}}, slightly modified]\label{thm:J_bounded}
Suppose that Assumptions~\ref{ass:J} and~\ref{ass:P_n_ker_J} are satisfied and the regularisation parameter $\alpha = \alpha(\delta,n)$ is chosen such that
	\begin{equation*}
	\alpha \to 0 \qand \frac{\left(\delta + \norm{A(I-P_{\U_n}) \Jminsol}\right)^2}{\alpha} \to 0 \quad \text{as $\delta \to 0$ and $n \to \infty$}.
	\end{equation*}
Then every sequence of minimisers $\unJ$ of ~\eqref{eq:proj_var_noisy} has a weakly convergent subsequence
\begin{equation*}
   \unJ \weakto \Jminsol.
\end{equation*}
\end{theorem}

Convergence rates can also be obtained under additional assumptions. 
Such rates are usually stated in terms of a (generalised) Bregman distance induced by the regularisation functional~\cite{Benning_Burger_modern:2018}.
We recall the definition for readers' convenience.
\begin{definition}[generalised Bregman distance]
For a proper convex functional $\reg$ the generalised Bregman distance between $x',x \in \U$ corresponding to the subgradient $q \in \dJ(x)$ is defined as follows
	\begin{equation*}
	D_\reg^q(x',x) \defeq \reg(x') - \reg(x) - \sp{q,x'-x}.
	\end{equation*}
	Here $\dJ(x)$ denotes the subdifferential of $\reg$ at $x \in \U$.
\end{definition}
The additional assumption required for obtaining a convergence rate is the \textit{source condition}.
\begin{theorem}[{\cite[Thm. 23]{asp-kor-sch-2020}}, slightly modified]\label{thm:conv_rates_var}
	Suppose that Assumptions~\ref{ass:J},~\ref{ass:P_n_ker_J} are satisfied and that $\Jminsol$ satisfies a source condition, i.e. there exists an element $q^\dagger \in \calY$ such that
	\begin{equation*}
	   A^*q^\dagger \in \dJ(\Jminsol).
	\end{equation*} 
	Then the following estimate holds for the Bregman distance between $\unJ$ and $\Jminsol$ 
	\begin{eqnarray*}
	D_\reg^{A^* q^\dagger}(\unJ,\Jminsol) &\leq&  \frac{1}{2\alpha} \left(\delta + \norm{A(I-P_{\U_n})\Jminsol} \right)^2 \\
    &+& \frac{\alpha}{2}\norm{q^\dagger}^2  +   \left(\delta \norm{q^\dagger} + C \norm{(I-P_{\U_n})A^*q^\dagger}\right)
	\end{eqnarray*}
	for some constant $C>0$. 
	
	If the regularisation parameter $\alpha = \alpha(\delta,n)$ is chosen as in Theorem~\ref{thm:J_bounded} then 
	\begin{equation*}
	D_\reg^{A^* q^\dagger}(\unJ,\Jminsol) \to 0 \quad \text{as $\delta \to 0$ and $n \to \infty$}.
	\end{equation*}
	For the particular choice
	\begin{equation}\label{var:parameter_choice}
	\alpha \sim \left(\delta + \max\left\{\norm{A(I-P_{\U_n})\Jminsol},\norm{(I-P_{\U_n})A^*q^\dagger}\right\}\right)
	\end{equation}
	we obtain the following estimate
	\begin{equation*}
	D_\reg^{A^* q^\dagger}(\unJ,\Jminsol)  \sim \alpha.
	\end{equation*}
\end{theorem}
The proofs of both theorems are similar to \cite{NeuSch90, Poeschl:2010}, although in those papers minimisers of the functional $x \in \calX \to \norm{Ax-y^\delta}^2+\alpha \reg(x)$ are approximated by  minimisers of the same functional over $\U_n$, while we solve the problem on the whole infinite-dimensional space.

We note that the convergence rate in \cref{thm:conv_rates_var} depends not only on how well the training images $\trim$ approximate the $\reg$-minimising solution $\Jminsol$, which is not surprising, but also on how well they approximate the subgradient $A^* q^\dagger$ from the source condition. This is another instance where training data for the adjoint operator $A^*$ may be advantageous.

We would also like to emphasise the different roles that the amount of training data $n$ plays in regularisation by projection (\cref{sec:reg-proj}) and variational regularisation. In regularisation by projection the solution is a linear combination of $n$ elements of the training set and therefore the size of this set $n$ controls the model complexity. The number of parameters in this case is the same as the number of training points. Furthermore, this number has to be controlled by the level of noise in the measurement $y^\delta$. In variational regularisation the solution is infinite-dimensional. Therefore, in some sense, we are in an overparametrised regime where the number of parameters (degrees of freedom in the solution) is infinite while the number of training points is finite. The parameter $n$ controls the approximation quality of the forward operator and can be chosen independently of the amount of noise in $y^\delta$.

%%%%%%%%%%%%%%%%%%%%%%%%%%%
\subsubsection{Iterative reconstruction methods and the role of the adjoint}\label{sec:proj-note-numerics}

The method~\eqref{eq:proj_var_noisy} has demonstrated good numerical performance learning and inverting the Radon transform~\cite[Sec. 6.4]{asp-kor-sch-2020}. These experiments used conic solvers provided by the CVX package~\cite{cvx}. Due to memory requirements, such solvers are not suited for large-scale applications, and iterative solvers are used instead. In this section we briefly discuss the application of such methods.

Perhaps the simplest iterative method, gradient descent, consists in taking the following updates for solving~\eqref{eq:proj_var_noisy}
\begin{multline}\label{eq:var-reg-GD}
    \iter{x}{k+1} = \iter{x}{k} - \tau_k ((AP_{\calX_n})^*(AP_{\calX_n})\iter{x}{k} - (AP_{\calX_n})^*y^\delta + \alpha q_k)  \\
            = \iter{x}{k} - \tau_k (P_{\calX_n}(A^*AP_{\calX_n}\iter{x}{k} - A^*y^\delta) + \alpha q_k), \quad q_k \in \dJ(\iter{x}{k}), 
\end{multline}
where  $\dJ(\iter{x}{k})$ is the subdifferential of $\reg$ at the iterate $\iter{x}{k}$ and $\tau_k > 0$ is the step size. 
We see that the iteration requires computing the operator $(AP_{\calX_n})^* = P_{\calX_n} A^*$. It is an easy calculation to show that it can be evaluated without numerical access to $A^*$:
\begin{equation}\label{eq:adj-APn}
    P_{\calX_n} A^* y = (AP_{\calX_n})^* y = \sum_{i=1}^n \sp{y, \orthim y^i} \orthim x^i, \quad y \in \calY.
\end{equation}

Compare~\eqref{eq:var-reg-GD} to the gradient descent step for the corresponding problem with the exact operator $A$,
\begin{equation}\label{eq:var-reg-GD-exact}
    \iter{x}{k+1} = \iter{x}{k} - \tau_k (A^* A\iter{x}{k} - A^*y^\delta + \alpha q_k).
\end{equation}

In~\eqref{eq:var-reg-GD} a projection $P_{\calX_n}$ is applied after the action of the restriction of the normal operator $A^*A P_{\calX_n}$. Depending on the problem, this may be a curse or a blessing. The range of the operator $A^*AP_{\calX_n}$ is the span of the training backprojections $\trback$,
\begin{equation*}
    \range{A^*AP_{\calX_n}} = \calZ_n = \span \trback,
\end{equation*}
while the projection $P_{\calX_n}$ will force the updates into the span of the training images $\trim$,
\begin{equation*}
    \range{P_{\calX_n} A^*AP_{\calX_n}} \subseteq \calX_n = \span \trim.
\end{equation*}
Depending on which subspaces, $\calX_n$ or $\calZ_n$, can better approximate the $\reg$-minimising solution $\Jminsol$, the projection $P_{\calX_n}$ in~\eqref{eq:var-reg-GD} may or may not be beneficial. We also note that if the forward operator is smoothing, then elements of $\calZ_n$ will be smoother than those of $\calX_n$.

The (outer) projection $P_{\calX_n}$ in~\eqref{eq:var-reg-GD} can be avoided if we have access to training data for the normal operator $\trback$, see~\eqref{eq:training-data}. In this case we can directly approximate the normal operator in~\eqref{eq:var-reg-GD-exact} with its restriction to $\calX_n$ and obtain
\begin{equation}\label{eq:var-reg-GD-normal}
    \iter{x}{k+1} = \iter{x}{k} - \tau_k (A^* A P_{\calX_n} \iter{x}{k} -  A^*P_{\calY_n} y^\delta + \alpha q_k).
\end{equation}
The operator $A^* A P_{\calX_n}$ can be evaluated without numerical access to either $A$ or $A^*$ using training data~\eqref{eq:training-data}. Indeed, applying $A^*$ to both sides of~\eqref{eq:APnx}, we get
\begin{equation*}
    A^* A P_{\calX_n} x = \sum_{i=1}^n \sp{x, \orthim x^i} \orthim z^i,
\end{equation*}
where $\trimhat$ are the orthonormalised training images and $\trbackhat$ the corresponding transformed backprojections. We need to apply a projection $P_{\calY_n}$ to the measurement $y^\delta$ before applying $A^*$ to match the range of the learned normal operator, $\range{A^*AP_{\calX_n}} = \calZ_n = \range{A^*P_{\calY_n}}$. The restricted operator $A^*P_{\calY_n}$ can also be evaluated without numerical access to $A^*$ using training data~\eqref{eq:training-data}. 

We will present numerical experiments with these approaches in \cref{sec:exp-proj-var}.

%%%%%%%%%%%%%%%%%%%%%%%%%%%%%%%%%%%%%%%%%%%%%%%%
%%%%%%%%%%%%%%%%%%%%%%%%%%%%%%%%%%%%%%%%%%%%%%%%
%%%%%%%%%%%%%%%%%%%%%%%%%%%%%%%%%%%%%%%%%%%%%%%%
\section{Data-driven model correction}\label{sec:model-correction}
In this section we will consider the case when an expensive forward model can be approximated by a  computationally more efficient one. For instance, in applications where the forward model is given by the solution of a partial 
differential equation, model reduction techniques are often used to reduce 
computational cost, e.g., by reduced order models or coarser discretisations. When the accurate model is replaced by a reduced one, this will lead to approximation errors, which may corrupt the reconstructed image depending on the severity of the approximation error. In the following, we will discuss how such model errors can be corrected with data-driven methods and  used to solve the inverse problem in a variational setting.

We recall that we consider linear inverse problems where we denote by $x\in \calX$  the unknown quantity of interest that we aim to reconstruct from measurements $y\in \calY$ 
and $x$ and $y$ fulfil the relation
\begin{equation}\label{eqn:invProb}
Ax=y,
\end{equation}
given bounded linear $A:\calX\to \calY$ (the accurate forward operator) acting between separable Hilbert spaces $\calX$ and $\calY$. 

We assume that the evaluation of the accurate operator $A$ is computationally expensive and we 
rather want to use a cheaper approximate model $\ApproxOp:\calX\to \calY$ with
\begin{equation}\label{eqn:invProb_approxModel}
    \ApproxOp x=\widetilde{y}, %=y+\delta y,
\end{equation}
leading to a systematic model error 
$
    \varepsilon = y - \widetilde{y}.
    $
In the following we will discuss different approaches to taking this systematic approximation error into account. First, we discuss a classical statistical correction, introduced as the approximation error method in~\cite{arridge:2006,kaipio2006statistical} and used in a specific application for linear corrections to nonlinear models in the variational setting \cite{arjas2023sequential}. We will then discuss  two  approaches, implicit and explicit corrections, in the framework of learned image reconstruction and the possibility to establish convergence guarantees for the explicit case \cite{lunz21}. 

\subsection{The \ac{AEM}} 
The well-established Bayesian \acl{AEM} \cite{kaipio2006statistical, arridge:2006} is an early data-driven approach to estimating a statistical model error. 
Recall that in Bayesian inversion we want to determine the posterior distribution of the unknown $x$ given $y$ and using Bayes' formula we obtain 
 \begin{equation}\label{eqn:Bayes}
 		p(x|y) = p(y|x) \frac{p(x)}{p(y)}.
 \end{equation}
 Thus, the posterior distribution is characterised by the likelihood $p(y|x)$ and the chosen prior $p(x)$ on the unknown. Typically, the likelihood $p(y|x)$ is modelled using accurate knowledge of the forward operator $A: \calX \to \calY$ as well as the noise model on the data $y$. The underlying idea of the approximation error method is  to adjust the likelihood by examining the difference between the (accurate) forward operator $A$ and its approximation $\ApproxOp$  \eqref{eqn:invProb}--\eqref{eqn:invProb_approxModel} as
\begin{equation}\label{eqn:modelErrorBAE}
    \varepsilon  =  Ax-\ApproxOp x. 
\end{equation}
Including an additive model for the measurement noise $e$, this leads to the modified observation model  
\begin{equation}\label{eqn:observationModelBAE}
    y =  \ApproxOp x + \varepsilon + e.
\end{equation} 
Here, both errors are (usually) assumed to be Gaussian. That is, first we model the measurement noise $e$ independent of $x$ as $e \sim \mathcal{N}(\eta_e,\Gamma_e)$,
where $\eta_e$ and $\Gamma_e$ are the mean and covariance.
Further, the model error $\varepsilon$ is approximated as Gaussian $\varepsilon \sim \mathcal{N}(\eta_{\varepsilon},\Gamma_{\varepsilon})$ and is assumed independent of the noise $e$ and the unknown $x$  leading to a Gaussian distributed total error $\totalerror=\varepsilon+e$, $\totalerror \sim \mathcal{N}(\eta_\totalerror,\Gamma_\totalerror)$,  where $\eta_{\varepsilon}$ and $\eta_\totalerror$ are means and $\Gamma_{\varepsilon}$ and $\Gamma_\totalerror$ are the covariance matrices of model error and total errors, respectively. We note here that the assumption of independence is a simplification and in practice one often observes dependence of the error on the signal.
Combining these leads to the so-called enhanced error model for the inverse problem \cite{kaipio2006statistical} with a likelihood distribution of the form
\[
    p(y|x) \sim \exp \left( - \frac{1}{2}\| L_\totalerror(\ApproxOp x - y + \eta_\totalerror) \|_\calY^2 \right)
\]
where $L_\totalerror$ such that $L_\totalerror^{\rm T}L_\totalerror=\Gamma_\totalerror^{-1}$  is a matrix square root such as the Cholesky decomposition of the inverse covariance matrix of the total error.
If the measurement noise $e$ is Gaussian white noise with  zero mean and constant standard deviation $\sigma$, the above formula can be written as
\[
p(y|x) \sim \exp \left( - \frac{1}{2\sigma}\| L_{\varepsilon}(\ApproxOp x - y + \eta_{\varepsilon}) \|_\calY^2 \right)
\]
where $L_{\varepsilon}^{\rm T}L_{\varepsilon}=\Gamma_{\varepsilon}^{-1}$. 
This motivates writing the variational problem, or the maximum a posteriori (MAP) estimator, for \eqref{eqn:Bayes} in the form
\begin{equation}\label{eqn:varProbBAE}
x^* = \argmin_{x\in X} \frac{1}{2}\|L_{\varepsilon} (\ApproxOp x - y + \eta_{\varepsilon}) \|_{\calY}^2 + \lambda \RegularisationOp(x),
\end{equation}
where $\RegularisationOp(x)$ is the regularisation functional corresponding to the prior $p(x)$ and $\lambda > 0$ is the regularisation parameter. 

In order to compute solutions, the unknown distribution of the model error needs to be approximated. 
This can be obtained for example by simulations \cite{arridge:2006,tarvainen2013bayesian} as follows. 
Let $\trim$ be the training set.  
The corresponding values of the model error are 
\begin{equation}
    \trni{\varepsilon}{i} = A\trni{x}{i}-\ApproxOp\trni{x}{i} 
\end{equation}
and the mean and covariance of the model error can  be estimated from the samples as 
\begin{eqnarray}
    \eta_{\varepsilon}  = \frac{1}{n} \sum_{i = 1}^{n} \trni{\varepsilon}{i} \, \text{ and }\, 
    \Gamma_{\varepsilon}  = \frac{1}{n-1} \sum_{i = 1}^{n} \trni{\varepsilon}{i}\otimes \trni{\varepsilon}{i} -\eta_{\varepsilon} \otimes \eta_{\varepsilon},
\end{eqnarray}
where $\otimes$ denotes the tensor product.

The approximation error method has found many applications in inverse problems, partly due to its simplicity yet high effectiveness in compensating for model errors. Examples of situations where it has been used successfully include model mismatch \cite{tarvainen2009approximation}, uncertainty in sensor locations \cite{sahlstrom2020modeling}, compensating for unknown boundary shapes \cite{candiani2021approximation}, and even recent applications in wireless communication \cite{marata2023joint}.

Despite these successes, it has been recently observed that the assumption of Gaussian distributed model errors as well as the independence assumption can be too restrictive, especially in nonlinear inverse problems \cite{smyl2021learning,arjas2023sequential,nicholson2023global}. This motivated the recent development of data-driven approaches for estimating non-Gaussian and non-independently distributed model errors as discussed in the following sections.

\subsubsection{Sequential model correction}
Let us first discuss briefly the work \cite{arjas2023sequential}, which further examines the non-Gaussianity of model errors in the case where the accurate forward model $A$ is nonlinear and the approximation $\ApproxOp$ is given by a linear model. This leads to a successive linearised and convexified problem that can be solved in a sequential manner as we will outline below.

In this case, we write \eqref{eqn:invProb} again in terms of $\ApproxOp$, which yields the observation model
\begin{equation} \label{eq:inverseproblemerror}
    y = \ApproxOp x + A(x) - \Tilde{A}x + e = \Tilde{A}x + \varepsilon(x) + e.
\end{equation}
We note that here the approximation creates a nonlinear approximation error, denoted by $\varepsilon(x)$, in contrast to the constant (i.e., independent of $x$) error in \eqref{eqn:observationModelBAE}. Consequently, this formulation of the model is still nonlinear as we have just moved the nonlinearity into $\varepsilon(x)$. Let us now assume that we have access to some initial reconstruction $x_0 \in \calX$. We can then approximate the model error by its value at $x_0$
\begin{equation}
    y \approx \Tilde{A}x + \varepsilon(x_0) + e.
\end{equation}
This leads to a convex variational formulation, given a convex regulariser $\RegularisationOp$,
\begin{equation} \label{eq:linprob}
    x^* = \argmin_{x \in \calX}\left\{\|\Tilde{A}x - (y - \varepsilon(x_0))\|_\calY + \lambda \RegularisationOp(x)\right\},
\end{equation}
which provides a local reconstruction depending on $x_0$. From here it is natural to expand this construction into a sequence 
\begin{equation} \label{eq:sequence}
    \iter{x}{k+1} = S(\iter{x}{k}) = \argmin_{x \in \calX}\left\{\|\Tilde{A}x - (y - \varepsilon(\iter{x}{k}))\|_\calY + \lambda \RegularisationOp(x)\right\}.
\end{equation}
We emphasise here that updating the sequence, i.e., solving the linearised and thus convex optimisation problem can be done efficiently with first-order optimisation methods. 

This approach is useful if a computationally effective linear model $\ApproxOp$ is available which makes \eqref{eq:sequence} tractable. The update for the sequence then only requires one evaluation of the the accurate nonlinear forward model and the solution of the linearised problem, which is cheaper then computing the Fr\'echet derivative of the nonlinear model $A$. The publication \cite{arjas2023sequential} shows that a fixed linear approximation $\ApproxOp$ performs well when compared to a scheme that updates the approximations between each sequential step. We also note that if one does not have access to the accurate model, or its evaluation even once is too expensive, a successive estimation of the model error could be computed or estimated \cite{smyl2021learning}. Such a sequential update of the approximation error is left for future research.

\subsection{Learned image reconstructions and implicit model corrections}
Let us now move to data-driven approaches in the context of learned image reconstruction. Here, broadly speaking, we aim to find a parameterised reconstruction operator $\mathcal{R}_\theta:\calY\to\calX$ whose parameters are learned from a suitable training set. This is most often achieved by utilising neural networks to parameterise the reconstruction operator; we refer to \cite{arridge_et_al_acta_numerica} for an overview of relevant methods. 

In what follows, we are interested in the framework of learned iterative reconstructions \cite{Adler_2017,hauptmann2018model}. That is, we aim to find a network $\Lambda_\Theta$ which is designed to mimic a gradient descent scheme. In particular, we train the network to perform an iterative update of the following form
\begin{equation}\label{eqn:LGS}
\iter{x}{k+1}=\Lambda_\Theta \left(\nabla_x \frac{1}{2}\|A\iter{x}{k}-y\|_{\calY}^2,\iter{x}{k}\right),
\end{equation}
where $\nabla_x \frac{1}{2}\|A\iter{x}{k}-y\|_{\calY}^2 = A^*(A\iter{x}{k}-y)$. 
If the accurate model is expensive to evaluate, computing the updates in \eqref{eqn:LGS} is expensive, which is especially problematic when training the networks. If the model is included in the training this quickly becomes intractable. 
Thus, one could use an approximate model $\ApproxOp$ instead of the accurate model and compute an approximate gradient as $\ApproxOp^*(\ApproxOp \iter{x}{k}-y)$ for the update in \eqref{eqn:LGS}, as proposed in \cite{hauptmann2018approximate}. The network $\Lambda_\Theta$ is then expected to \textit{implicitly} correct the introduced model error to produce a new reliable iterate.

That means that correction and regularisation are  trained simultaneously with the update in \eqref{eqn:LGS}. Such approaches are typically trained by using a loss function, such as the $L^2$-loss, to measure the distance between reconstruction and a ground truth phantom. This way a substantial speed-up, compared to classical variational approaches, can be achieved with improved reconstruction quality. In a recent paper \cite{hauptmann2023model} the implicit correction has been extended to a model corrected learned primal dual method, where separate updating operators are learned in primal and dual space, offering further improvements of reconstruction quality.

Nevertheless, such implicit corrections within a learned reconstruction operator offer limited insights into how approximate models are corrected for and so far have only limited convergence guarantees \cite{mukherjee2023learned}. Thus, we will consider in the following an \textit{explicit} correction that can  be subsequently used in a variational framework.

\subsection{Explicit model correction and a convergence result}\label{sec:explicit}
Let us now consider corrections for the approximation error caused by the approximate model $\ApproxOp$ via a 
parameterisable nonlinear mapping $\ForwardCor: \calY\to \calY$, applied directly as correction to $\ApproxOp$ as proposed in \cite{lunz21}. This mapping could be given by a (convolutional) neural network, but other options can be considered. This leads to a corrected operator $\CorrectedOp$ of the form
\begin{equation}\label{eqn:correctedModel}
    \CorrectedOp=\ForwardCor\circ\ApproxOp.
\end{equation}
We aim to choose the correction $\ForwardCor$ such that ideally $\CorrectedOp(x) \approx \TrueOp x$ for those $x\in \calX$ that we are interested in. 
The primary question that we aim to answer is, whether such corrected 
models \eqref{eqn:correctedModel} can be subsequently used in variational regularisation
approaches. 
Thus, it is natural to require that the obtained solutions involving the corrected operator $\CorrectedOp$ and the accurate operator $\TrueOp$, are close, that is 
\begin{equation}\label{equ:variationalSolutionEquality}
    \argmin_{x\in \calX} \frac{1}{2} \|\CorrectedOp(x) - y\|_{\calY}^2 + \lambda \RegularisationOp(x)
    \approx
    \argmin_{x\in \calX}  \frac{1}{2} \|\TrueOp x - y\|_{\calY}^2 + \lambda \RegularisationOp(x).
\end{equation} 
Solutions are then usually computed by an iterative algorithm.
Here we consider first order methods to draw connections to learned iterative schemes as in \eqref{eqn:LGS}. In particular, we consider a classical gradient descent scheme, assuming a differentiable $\RegularisationOp$. Given an 
initial guess $x_0$, we  compute a solution by the iterative process
\begin{equation}\label{equ:proximalGrad}
    \iter{x}{k+1} = \iter{x}{k} - \gamma_k \nabla_x\left(\frac{1}{2}\|\CorrectedOp \iter{x}{k}-y\|_{\calY}^2 +\lambda \RegularisationOp(\iter{x}{k})\right),
\end{equation}
with an appropriately chosen step size $\gamma_k>0$. When using \eqref{equ:proximalGrad} for the corrected operator it seems natural to ask for
a \textit{gradient consistency} {of the approximate gradient}
$
    \nabla_x \|\CorrectedOp(x) - y\|_{\calY}^2 \approx \nabla_x \|\TrueOp x - y\|_{\calY}^2.
$
We recall that the correction $\ForwardCor$ in \eqref{eqn:correctedModel} is given by a nonlinear neural network and with the chain rule we obtain 
\begin{equation}\label{eqn:fidelityTerm}
     \frac{1}{2}\nabla_x \|\CorrectedOp(x) - y\|_\calY^2 = \ApproxOp^*  \left[D\ForwardCor(\ApproxOp x)\right]^* \left(\ForwardCor(\ApproxOp x)-y\right).
\end{equation}
Here, we denote by $D\ForwardCor(y)$ the Fr{\'{e}}chet derivative of $F_\Theta$ at $y$, which is a linear operator $\calY \to \calY$. 
That means, to satisfy the gradient consistency condition, we would need
\begin{equation}\label{equ:gradConsistencyApprox}
\ApproxOp^* \left[D\ForwardCor(\ApproxOp x)\right]^* \left(\ForwardCor (\ApproxOp x)-y\right) \approx A^*(Ax-y).
\end{equation}
This reveals a problem: the range of the corrected fidelity term's gradient \eqref{eqn:fidelityTerm} is limited  by the range of the approximate 
adjoint, $\range{\ApproxOp^*}$. Thus, we identify the key difficulty here in the differences of the range of the accurate and
the approximate adjoints rather than the differences in the forward operators themselves.
A correction of the forward operator via composition with a parametrised model $F_\Theta$ in measurement space is not able to produce gradients close to the gradients of the accurate data term if $\range{\ApproxOp^*}$ and $\range{A^*}$ are too different, see also Theorem 3.1 in \cite{lunz21}.

\subsubsection{Obtaining a Forward-Adjoint Correction}\label{sec:forw-adj-cor}
To achieve a gradient--consistent model correction two networks can be considered instead. That is, we learn a network $\ForwardCor$ that 
corrects the forward model and another network $\AdjointCor$ that corrects the adjoint, such that we have
\begin{equation}
    \CorrectedOp := \ForwardCor \circ \ApproxOp, \quad \CorrectedAd := \AdjointCor \circ \ApproxOp^*
\end{equation}
These corrections can then be obtained as follows. Given a set of training samples $\{\trni{x}{i}, \trni{y}{i} = A\trni{x}{i}\}_\ind{i}{1}{n}$, we train the forward correction $\ForwardCor$ acting in measurement space $\calY$, for the adjoint we train the network $\AdjointCor$
acting on image space $\calX$ using  two losses 
\begin{equation}\label{equ:forwardLoss_FBC}
    \min_\Theta \sum_i \| \ForwardCor(\ApproxOp \trni{x}{i}) - A \trni{x}{i} \|_{\calY} \quad \text{and} \quad \min_\Phi \sum_i \| \AdjointCor(\ApproxOp^* \trni{r}{i}) -  A^* \trni{r}{i} \|_{\calX},
\end{equation}
where $\trni{r}{i} = \ForwardCor( \ApproxOp \trni{x}{i})- \trni{y}{i} $.
This ensures that the adjoint correction is in fact trained in directions relevant when solving the variational problem. We can then use both corrections to compute approximate gradients of the data fidelity term $\|Ax-y\|_{\calY}^2$ as
\begin{equation}
    \label{eqn:GradientFit}
    A^*(Ax-y)  \approx \left(\AdjointCor\circ \ApproxOp^*\right) \left( \ForwardCor(\ApproxOp x) - y \right).
\end{equation}
A convergence result can be established by considering the two functionals corresponding to accurate and corrected operator as
\begin{equation}
    \mathcal{L}(x) := \frac{1}{2} \| \TrueOp x - y \|_{\calY}^2 + \lambda \RegularisationOp(x), \, 
    \mathcal{L}_\Theta(x) := \frac{1}{2} \| \CorrectedOp(x) - y \|_{\calY}^2 + \lambda \RegularisationOp(x)
\end{equation}
and using the forward-adjoint correction in the minimisation. We can then obtain, under suitable conditions outlined in \cite{lunz21}, the following theorem.
\begin{theorem}[Convergence to a neighbourhood of the accurate solution \cite{lunz21}] \label{thm:conv-neighbourhood}
    Let $\epsilon>0$ and suitable constant $C>0$ (controlling the subdifferential of $\mathcal{L}_\Theta$). Assume both adjoint and forward operator are fit up to a $C/4$-margin, i.e.
    \begin{equation}\label{eqn:conditionTheo410}
       \|\TrueOp\|_{\calX \to \calY} \| (\TrueOp - \CorrectedOp)(\iter{x}{k}) \|_{\calY} < C/4, \quad \|(\TrueOp^* - \CorrectedAd) (\CorrectedOp(\iter{x}{k}) - y) \|_{\calX} < C/4
    \end{equation}
    for all $y$ and $\iter{x}{k}$ obtained during gradient descent over $\mathcal{L}_\Theta$. 
    Then eventually the gradient descent dynamics over $\mathcal{L}_\Theta$ will reach an $\epsilon$ neighbourhood of the solution of the problem corresponding to the exact operator.    
\end{theorem}

The proof of Theorem \ref{thm:conv-neighbourhood} relies on ensuring that the gradients \eqref{eqn:GradientFit} are pointing in the same direction to yield a descent direction with respect to the accurate functional. That is, we want to ensure that the angle between the approximate and the exact gradients is positive,
\begin{equation}
    \label{eqn:AligenementVariational}
    \cos \Phi_v(x) := \frac{\langle \nabla \mathcal{L}(x), \nabla^\dagger  \mathcal{L}_\Theta(x) \rangle}{\| \nabla \mathcal{L}(x)\|^2}>0,
\end{equation}
where $\nabla^\dagger$ is used to indicate that we compute a corrected gradient of $\mathcal{L}_\Theta(x)$ using the forward-adjoint correction given by the right-hand side of \eqref{eqn:GradientFit}. The experiments in \cite{lunz21} show that when this alignment is ensured during the minimisation then indeed one can observe  convergence to the same neighbourhood as with the accurate model, while if a positive alignment is not ensured then the optimisation procedure will diverge. 

In~\cite{lunz21} it was proposed to use the so-called \textit{recursive} training to ensure this positive alignment. It requires running the iterations for every training image $x^i$ and adding the values of the forward operator and its adjoint along the trajectory to the training set. This can be  very expensive computationally.

%%%%%%%%%%%%%%%%%%%%%%%%%%%%%%%%%%%%%%%%%%%%%%%%
\subsubsection{Limitations and extensions}\label{sec:GD}

The first, minor difficulty with this approach is that it requires training two networks, one for the forward operator $A$ and one for its adjoint $A^*$. This can be avoided by using training data for the normal operator $A^*A$ as in \cref{sec:proj-note-numerics} and learning a single network $N_\Theta \colon \calX \to \calX$ to satisfy $N_\Theta(\tilde A^* \tilde A x^i) \approx A^*Ax^i$. This will be sufficient to approximate the gradient of the data fidelity term $\|Ax-y\|_{\calY}^2$, see~\eqref{eqn:GradientFit}.

A much more serious difficulty is that \cref{thm:conv-neighbourhood} requires that the trained networks approximate the exact forward and adjoint operators for all iterates $\iter{x}{k}$ obtained during gradient descent (and not just in the vicinity of training samples $\trim$). This requires very cumbersome and computationally expensive recursive training discussed above.

A possible remedy for this is to use iterative algorithms that stay in the vicinity of the training images $\trim$, such as the projected gradient scheme in~\cite{hauswirth16} where the authors consider a constrained variational problem over a manifold. A difficulty is, however, that in our setting the manifold is given implicitly via training samples $\trim$ and needs to be estimated ``on the fly'' as the iterations proceed. This is currently work in progress.

%%%%%%%%%%%%
\subsubsection{Connections to regularisation by projection}

The simplified model $\tilde A$ can also be used in the context of regularisation by projection (see \cref{sec:proj-var-reg}). While the learned linear model $AP_{\calX_n}$ is exact on the subspace $\calX_n$ spanned by the training images $\trim$, on the complement of this subspace the learned model is zero, and it may be beneficial to use the simplified model $\tilde A$ instead. This will lead to the following approximation of the forward operator
\begin{equation*}
    A \approx AP_{\calX_n} + \tilde A(Id-P_{\calX_n})
\end{equation*}
The corresponding  variational problem will then read as follows (cf.~\eqref{eq:proj_var_noisy})
\begin{equation}\label{eq:proj_var_combined}
\min_{x \in \U} \frac12 \norm{\left[AP_{\U_n} + \tilde A (Id-P_{\U_n})\right] x - y^\delta   }^2 + \alpha \reg(x)
\end{equation}
and the gradient descent iteration will become
\begin{multline}\label{eq:var-reg-GD-combined}
    \iter{x}{k+1} = \iter{x}{k} - \tau_k \left( \left[P_{\calX_n} A^*AP_{\calX_n} + (Id-P_{\calX_n}) \tilde A^* \tilde A (Id-P_{\calX_n})\right] \iter{x}{k} \right.  \\
    - \left. \left(P_{\calX_n}A^* + (Id-P_{\calX_n})\tilde A^* \right) y^\delta + \alpha q_k\right), \quad q_k \in \dJ(\iter{x}{k}).
\end{multline}
All operators here can be computed without numerical access to the exact model $A$, relying only on the training pairs $\{x^i,y^i=Ax^i\}_{i=1,...,n}$ and the simplified model $\tilde A$. (Recall that the operator $P_{\calX_n} A^* = (AP_{\calX_n})^*$ can be computed using these training pairs via~\eqref{eq:adj-APn}).

Alternatively, the iteration~\eqref{eq:var-reg-GD-normal}, based on learning the normal operator on $\calX_n$, can be augmented with an approximate component on $\calX_n^\perp$
\begin{multline}\label{eq:var-reg-GD-normal-combined}
    \iter{x}{k+1} = \iter{x}{k} - \tau_k \left(\left[A^* A P_{\calX_n} + \tilde A^* \tilde A (Id-P_{\calX_n})\right]\iter{x}{k} \right. \\
    - \left. \left(A^*P_{\calY_n} + \tilde A^* (Id-P_{\calY_n}) \right)y^\delta + \alpha q_k \right).
\end{multline}
This iteration only requires the training data $\{x^i,z^i=A^*Ax^i\}_{i=1,...,n}$ for the normal operator and the simplified model $\tilde A$ along with its adjoint $\tilde A^*$.

Numerical experiments with these approaches will be presented in \cref{sec:exp-proj-var}.

%%%%%%%%%%%%%%%%%%%%%%%%%%%%%%%%%%%%%%%%%%%%%%%%
\section{\ac{PAT}}\label{sec:PAT}
In this section we will apply the presented approaches to \acf{PAT}. To do that, let us first briefly discuss the \ac{PAT} forward problem and then introduce the analytic approximate model for the context of model corrections.

To create the measured signal in PAT, biological tissue is exposed to a sufficiently short near-infrared light pulse that is then absorbed by chromophores. This results in a spatially-varying pressure increase, which initiates an 
\ac{US} pulse, that then propagates to the tissue surface. The measurement consequently consists of the detected waves in space-time at the boundary of the tissue. 
This time evolution of the  photoacoustic wave can be modelled using the equations of linear acoustics \cite{cox2005fast,treeby2010kWave}, and can be described as an initial value problem for the wave equation with spatial coordinates $\zeta\in\R^2$ and time $t\geq 0$  
\begin{subequations}
\begin{eqnarray}
    \label{eq:wave_equation_pat}
   (\partial_{tt} - c^2 \Delta)p(\zeta,t) &=& 0, \\
  \label{eq:wave_initialcondition1}
   p(\zeta,t=0)   &=& p_0(\zeta), \\
  \label{eq:wave_initialcondition2}
   \partial_t p(\zeta,t=0) &=& 0,
\end{eqnarray}
\end{subequations}
where $c$ is the speed of sound. The measurement of the time series is  modelled as a linear operator $\mathcal{M}$ acting on the pressure field $p(\zeta,t)$ restricted to the boundary $\Gamma$ of the computational domain and a finite time window:
\begin{equation}
y = \mathcal{M} \, p_{|\Gamma \times (0,T)}. \label{eqn:Measurement}
\end{equation}
Together, equations \eqref{eq:wave_equation_pat} and \eqref{eqn:Measurement} define the linear forward model $A$ that maps
the initial pressure $x=p_0(\zeta)$ to the measured time series $y$. This forward model can be accurately simulated by a pseudo-spectral time-stepping model as outlined in \cite{treeby2010kWave,treeby2012modeling}. While providing an efficient implementation, time-stepping can  take a considerable amount of time depending on a possibly fine time discretisation. 

\begin{figure}[t!]
    \centering
   \includegraphics[width=0.8\textwidth]{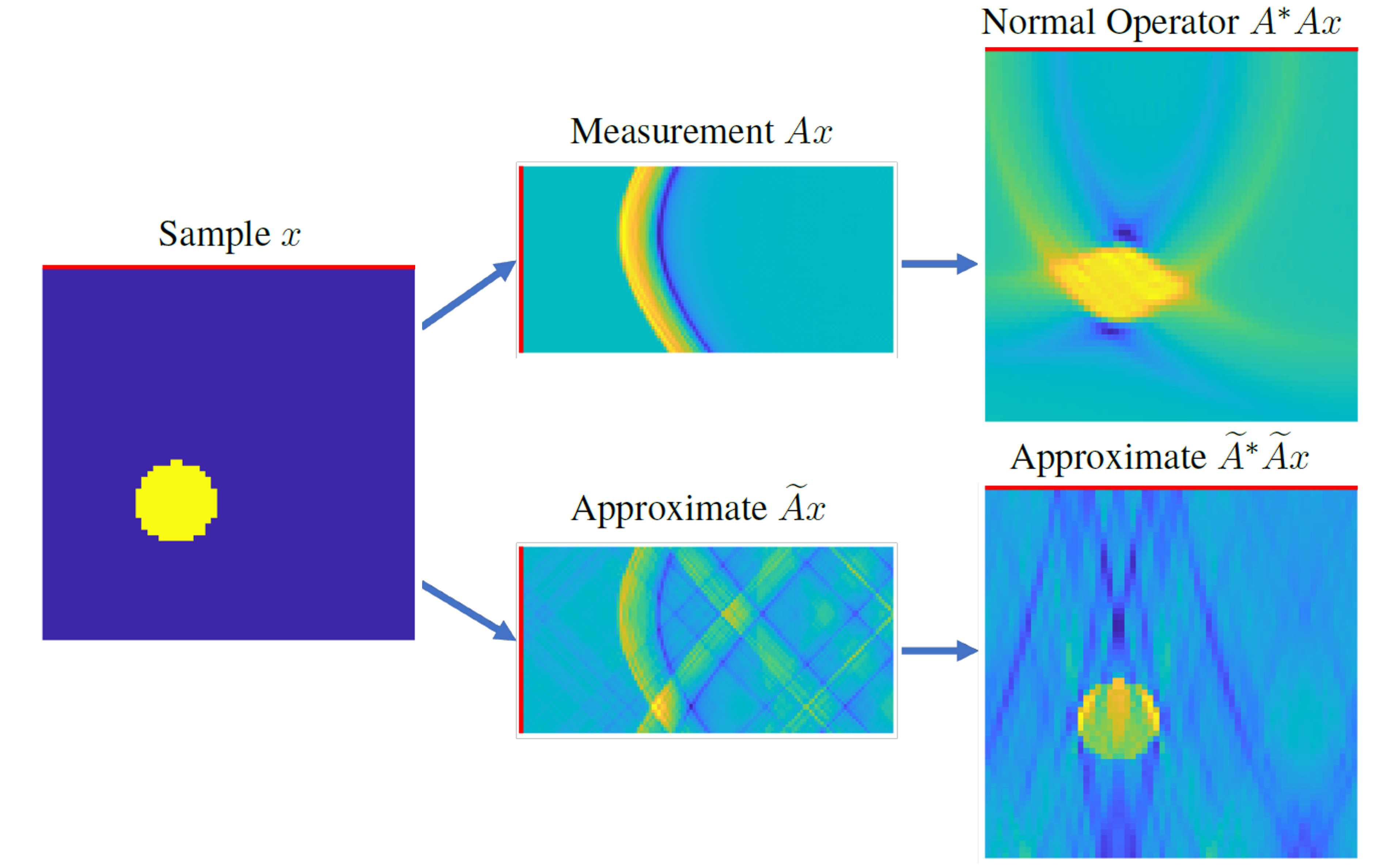}
 \caption{\label{fig:ForwApprox} 
 Illustration of the forward operator $A$ and the approximate model $\widetilde A$, as well as corresponding normal operators. The sensor is indicated with a red line (left/right images), in the measurements (middle) the red line corresponds to initial time $t=0$. }
\end{figure}

Thus, we can consider a model that eliminates the time stepping and replaces it with one (Fast) Fourier transform. First, we can consider the case where measurement points lie on a line ($\zeta_2=0$) outside the support of $x$. Then under the assumption of constant speed-of-sound, the pressure on the sensor can be related to $x$ as follows \cite{cox2005fast,kostli2001}
\begin{equation}
p(\zeta_1,t) = \frac{1}{c^2} 
\mathcal{F}^{-1}_{k_1} \left\{\mathcal{C}_{\omega}\left\{
B({k_1},\omega) \tilde{x}({k_1},\omega)
\right\}\right\}, 
\label{eqn:FastFwd}
\end{equation}
where 
$\tilde{x}({k_1},\omega)$ is obtained from the Fourier transform $\hat{x}(k) = \mathcal{F}_{\zeta}\{x(\zeta)\}$ via the dispersion relation $(\omega/c)^2 = k_1^2+k_2^2$,  $\mathcal{C}_{\omega}$ is a cosine transform from $\omega$ to $t$, and $\mathcal{F}^{-1}_{k_1}$ is the 1D inverse Fourier Transform from $k_1$ to $\zeta_1$ on the detector. The weighting factor,
\begin{equation}
B(k_1,\omega) = \omega/\left(\textrm{sgn}(\omega)\sqrt{(\omega/c)^2 - k_1^2}\right),
\end{equation}
contains an integrable singularity which means that if \eqref{eqn:FastFwd} is evaluated by discretisation on a rectangular grid   (enabling the application of FFT for efficient calculation), then aliasing will apear in the measured data $p(\zeta_1,t)$. Consequently, evaluating \eqref{eqn:FastFwd} using FFT leads to a \textit{fast but approximate} forward model. We can control the degree of aliasing by avoiding waves that arrive close to parallel at the sensor. This could be included in the model as an angular thresholding to control the degree of aliasing, we refer to \cite{hauptmann2018approximate} to a more detailed discussion. 
Either way, with or without angular thresholding in the weighting factor $B$, the relation \eqref{eqn:FastFwd} defines the approximate model $\tilde A$ used in what follows. The difference between the accurate model $A$ and the approximate model $\widetilde{A}$ is shown in Figure \ref{fig:ForwApprox}. The aliasing artefacts are clearly visible in the data space and the resulting normal operator does carry wrong information for the reconstruction.

%%%%%%%%%%%%%%%%%%%%%%%%%%%%%%%%%%%%%%%%%%%%%%%%
\subsection{Training data}\label{sec:training-data}
We will consider a computational domain
$(\zeta_1,\zeta_2) \in \Omega = [0, 1] \times [0, 1]$ with a rectangular discretisation of $64 \times 64$ pixels. The measurements
are taken  at the top  of the domain. 
The background value is set to zero and we sample a number of indicator functions of discs located randomly in the domain with
parameters uniformly distributed as follows: centre $(\zeta_1, \zeta_2) \in [0.25, 0.75] \times [0.25, 0.75]$, radius $r \in [0.1, 0.2]$. Each disc has a random  initial pressure value $p_0 \in [0.5, 1]$.

%%%%%%%%%%%%%%%%%%%%%%%%%%%%%%%%%%%%%%%%%%%%%%%%
\subsection{Experiments with projected variational regularisation}\label{sec:exp-proj-var}

In this section we present numerical experiments with the method described in \cref{sec:proj-var-reg}. 
As the regulariser $\reg$ we take the Total Variation ($\TV$), which we define as a functional on $L^2(\Omega)$ extending it with the value $+\infty$ on $L^2(\Omega) \setminus \BV(\Omega)$. 
This is a common setting in imaging~\cite{Chambolle_Lyons:1997}. 

To orthonormalise the training images $\trim$, we use the modified Gram-Schmidt algorithm~\cite{trefethen_num_lin_al}. Due to its numerical instability, we restrict ourselves to $n=1500$ training samples. 

\begin{figure}[th!]
    \captionsetup[subfigure]{justification=centering}		
    \begin{minipage}{0.8\textwidth}
    \centering
    ~~~~~~\begin{subfigure}[t]{0.4\textwidth}
    \begin{tikzpicture}
    \node[inner sep=0] (image) at (0,0) {\adjincludegraphics[width=\textwidth,trim={0 0 0 {.07\height}},clip]{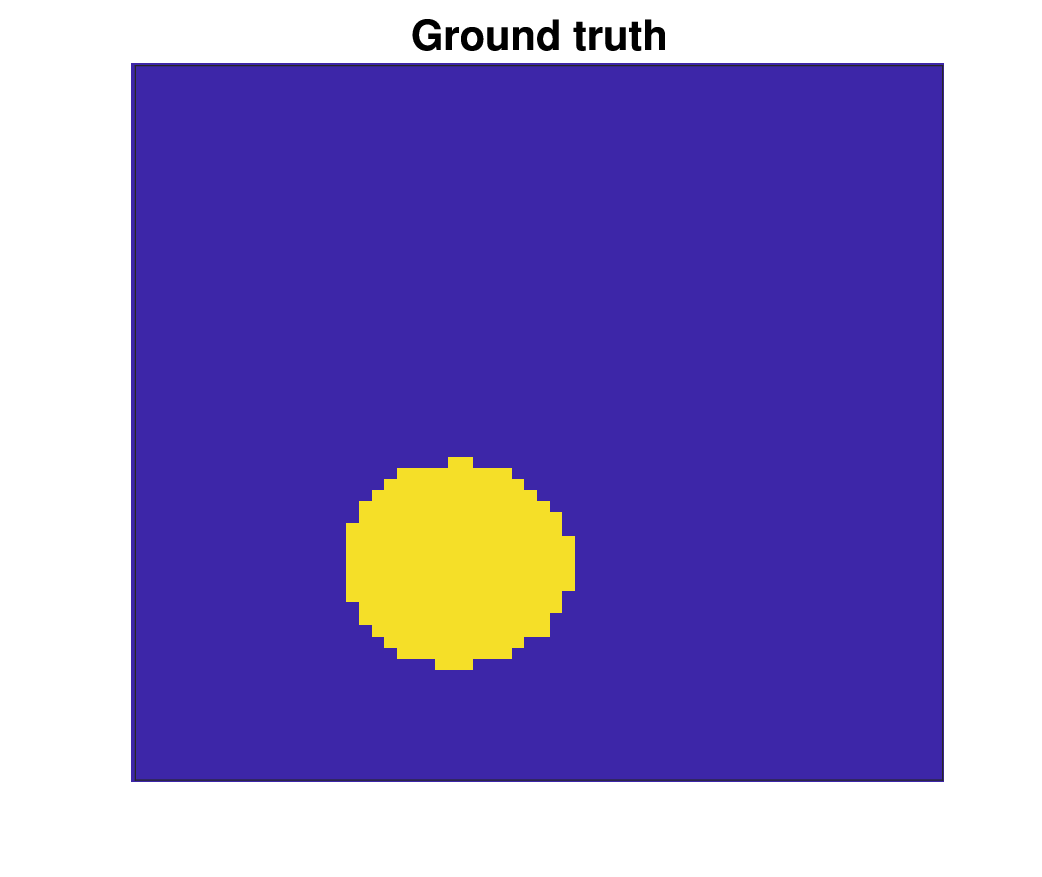}};
    \node[above=0 of image] {\scriptsize Ground truth};
    \end{tikzpicture}
	\end{subfigure}
	\begin{subfigure}[t]{0.4\textwidth}
    \begin{tikzpicture}
    \node[inner sep=0] (image) at (0,0) {\adjincludegraphics[width=\textwidth,trim={0 0 0 {.07\height}},clip]{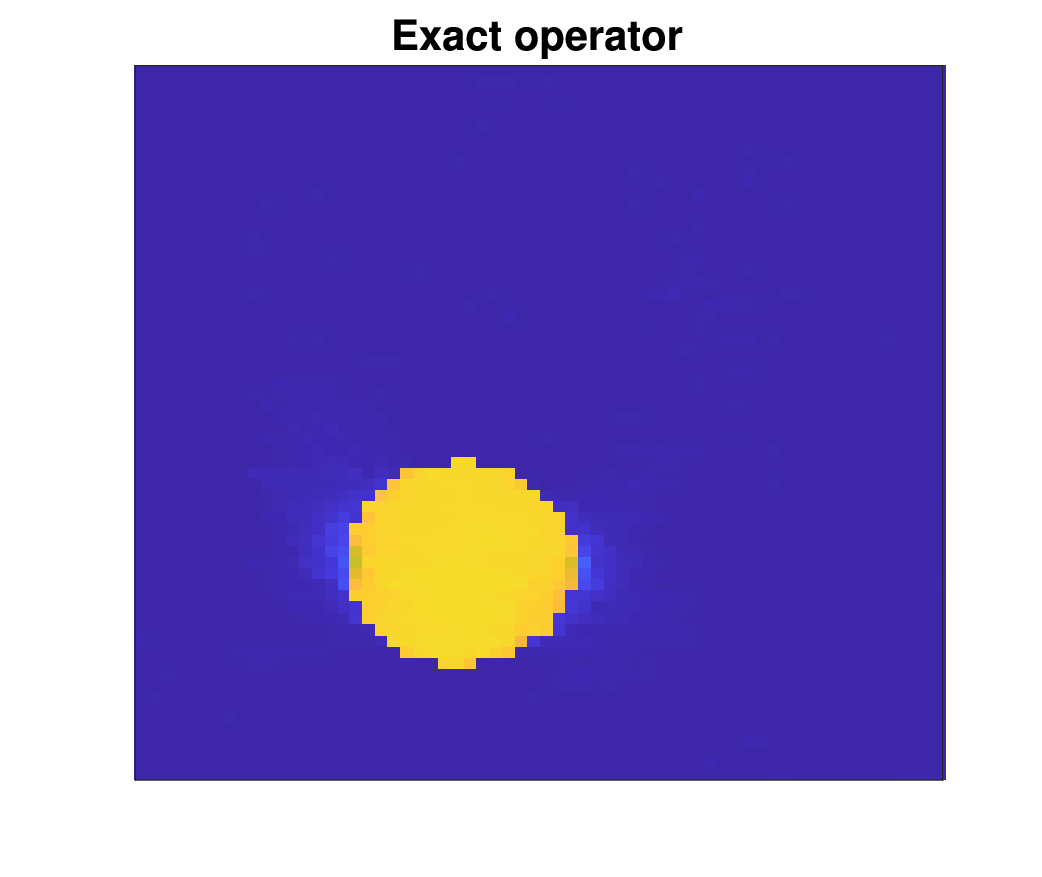}};
    \node[above=0 of image] {\scriptsize Exact operator};
    \end{tikzpicture}
    \end{subfigure} \\[5pt]
    ~~~~~~\begin{subfigure}[t]{0.4\textwidth}
    \begin{tikzpicture}
    \node[inner sep=0] (image) at (0,0) {\adjincludegraphics[width=\textwidth,trim={0 0 0 {.07\height}},clip]{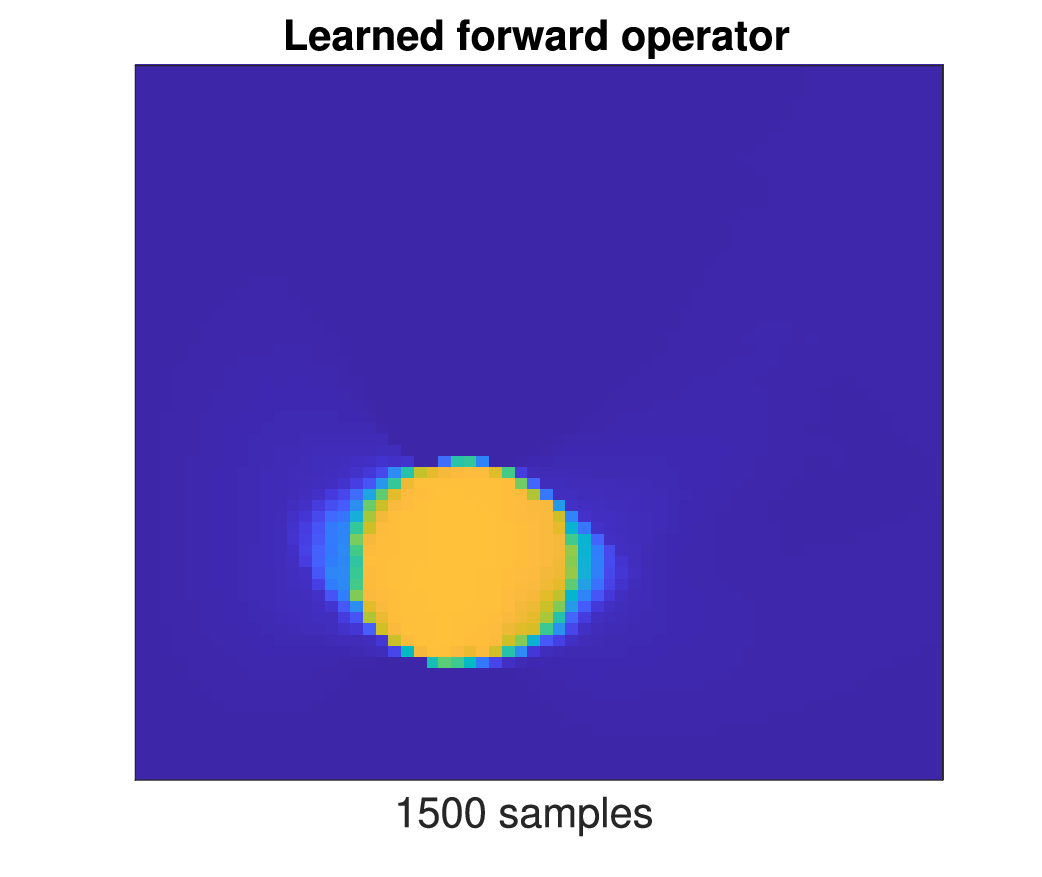}};
    \node[above=0 of image] {\scriptsize Learned forward operator};
    \end{tikzpicture}
	\end{subfigure}
	\begin{subfigure}[t]{0.4\textwidth}
    \begin{tikzpicture}
    \node[inner sep=0] (image) at (0,0) {\adjincludegraphics[width=\textwidth,trim={0 0 0 {.07\height}},clip]{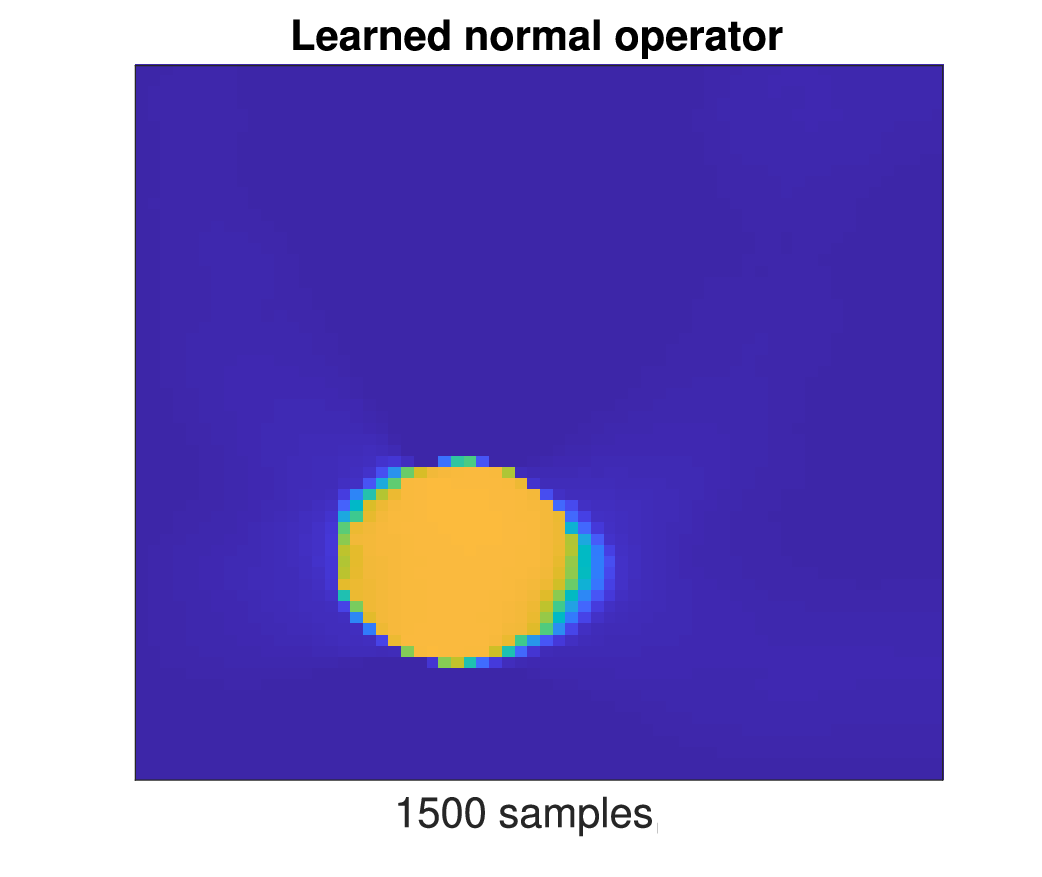}};
    \node[above=0 of image] {\scriptsize Learned normal operator};
    \end{tikzpicture}
    \end{subfigure} \\[5pt]
    ~~~~~~\begin{subfigure}[t]{0.4\textwidth}
    \begin{tikzpicture}
    \node[inner sep=0] (image) at (0,0) {\adjincludegraphics[width=\textwidth,trim={0 0 0 {.07\height}},clip]{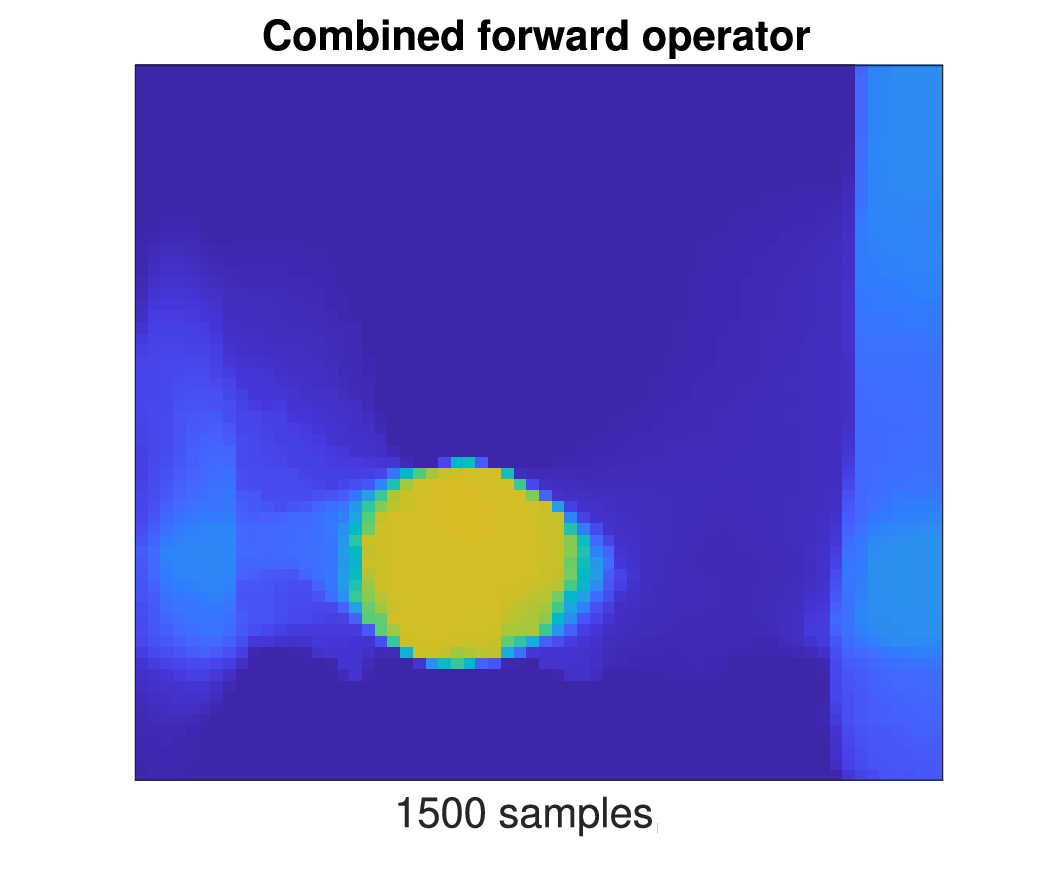}};
    \node[above=0 of image] {\scriptsize Combined forward operator};
    \end{tikzpicture}
	\end{subfigure}
	\begin{subfigure}[t]{0.4\textwidth}
    \begin{tikzpicture}
    \node[inner sep=0] (image) at (0,0) {\adjincludegraphics[width=\textwidth,trim={0 0 0 {.07\height}},clip]{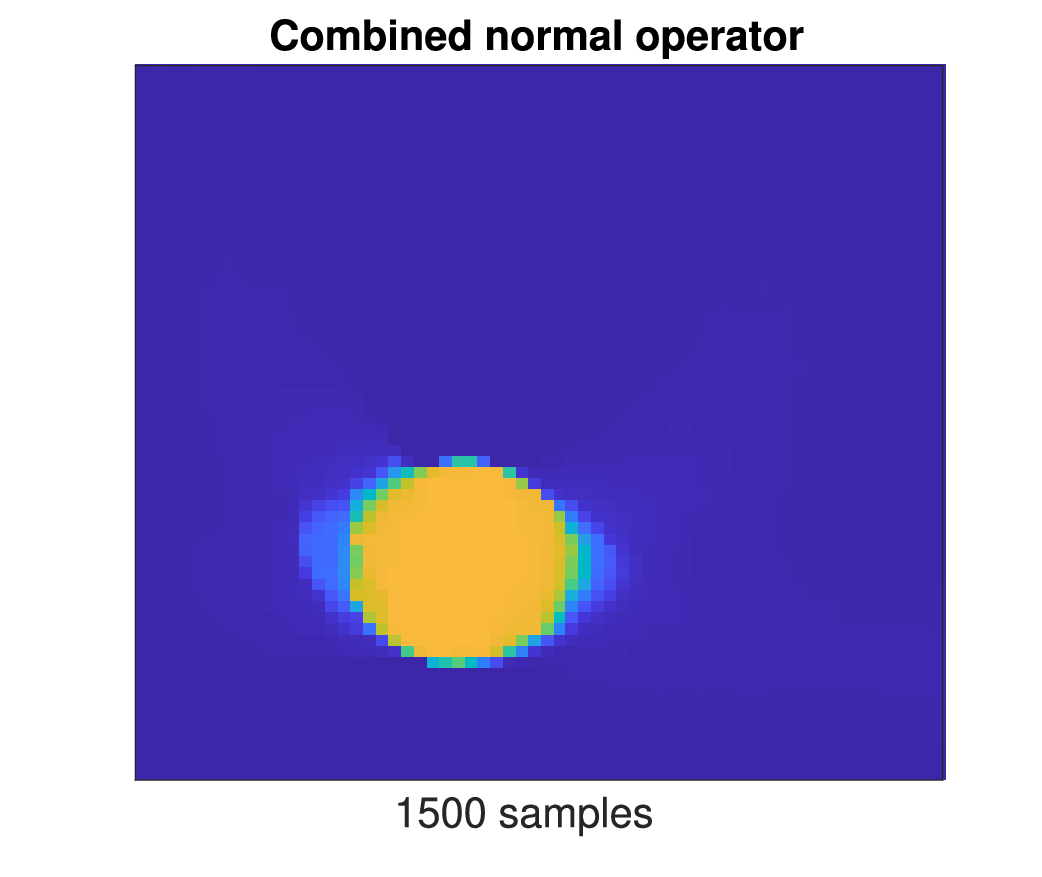}};
    \node[above=0 of image] {\scriptsize Combined normal operator};
    \end{tikzpicture}
    \end{subfigure}
    \end{minipage}%
    \begin{minipage}{0.2\textwidth}
        \includegraphics[width=0.5\textwidth]{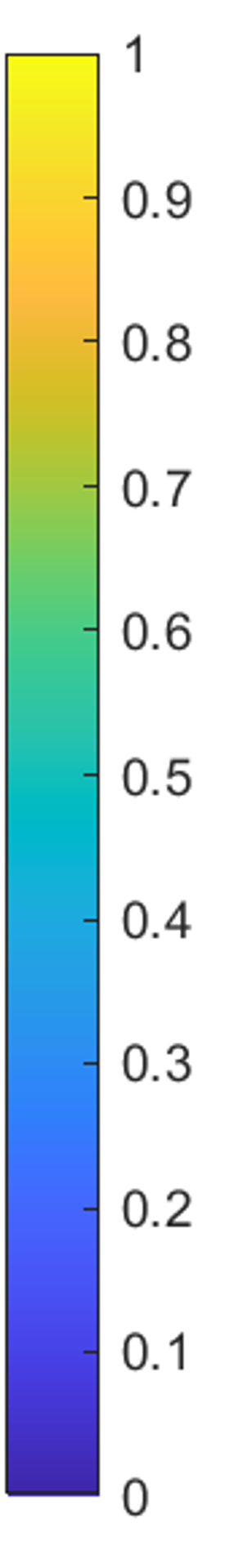}
    \end{minipage}
    \caption{Projected variational regularisation. Experiments with $n=1500$ training samples. Both learned forward operator and learned normal operator perform well, the reconstructions with the learned normal operator being perhaps a bit sharper. Surprisingly, combining the learned forward model with an approximate one $\tilde A$ decreases the reconstruction quality.}\label{fig:proj-var}
\end{figure}
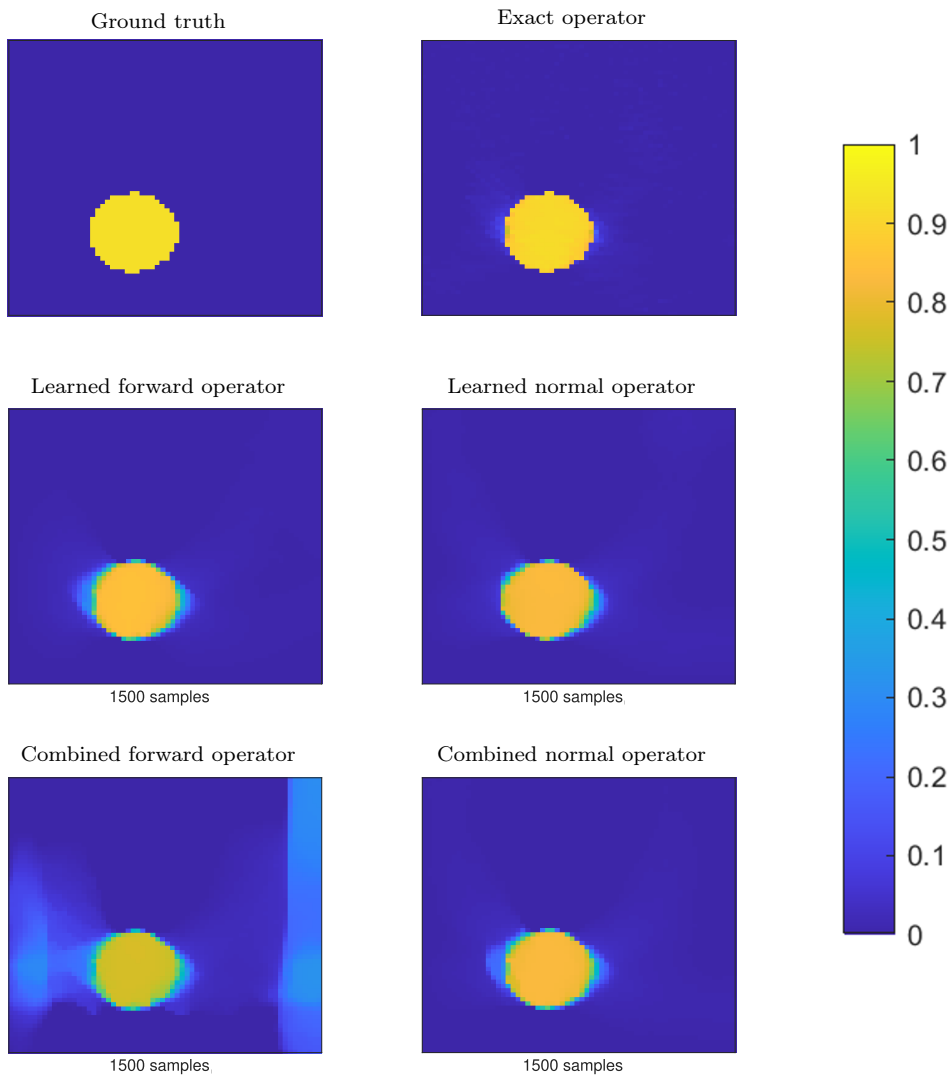

We use the iterations \cref{eq:var-reg-GD} -- learned forward operator, \eqref{eq:var-reg-GD-normal} -- learned normal operator, \cref{eq:var-reg-GD-combined} -- learned forward operator combined with approximate model, and~\eqref{eq:var-reg-GD-normal-combined} -- learned normal operator combined with approximate model. The results are shown in \cref{fig:proj-var}. The top row shows the ground truth and the reconstruction with the exact operator, which is our golden standard in these experiments. Performing a basic parameter search, we find a reasonable value of the regularisation parameter $\alpha=2 \cdot 10^{-4}$ which yields a relative reconstruction error of $6\%$.

With a learned forward (middle row left) and a learned normal operator (middle row right) we obtain a relative reconstruction error of $23-24\%$, also after performing a basic parameter search and finding a reasonable value of $\alpha=2 \cdot 10^{-2}$. The value of the parameter is higher, as expected in a problem with operator errors. The reconstruction obtained with a learned normal operator looks sharper than the one obtained with a learned forward operator, but overall the reconstruction quality is comparable. 

\begin{figure}[t!]
    \captionsetup[subfigure]{justification=centering}
    \begin{minipage}{0.8\textwidth}
    \centering
    ~~~~~~\begin{subfigure}[t]{0.4\textwidth}
    \begin{tikzpicture}
    \node[inner sep=0] (image) at (0,0) {\adjincludegraphics[width=\textwidth,trim={0 0 0 {.07\height}},clip]{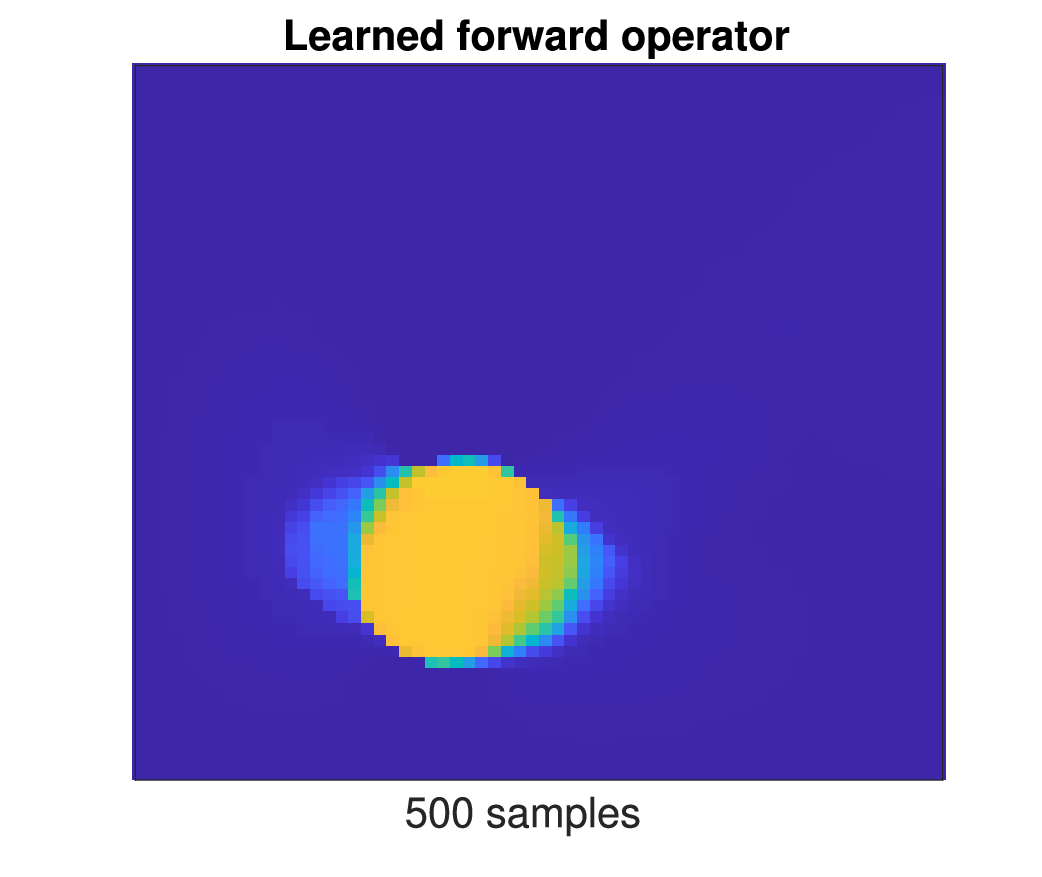}};
    \node[above=0 of image] {\scriptsize Learned forward operator};
    \end{tikzpicture}
	\end{subfigure}
	\begin{subfigure}[t]{0.4\textwidth}
    \begin{tikzpicture}
    \node[inner sep=0] (image) at (0,0) {\adjincludegraphics[width=\textwidth,trim={0 0 0 {.07\height}},clip]{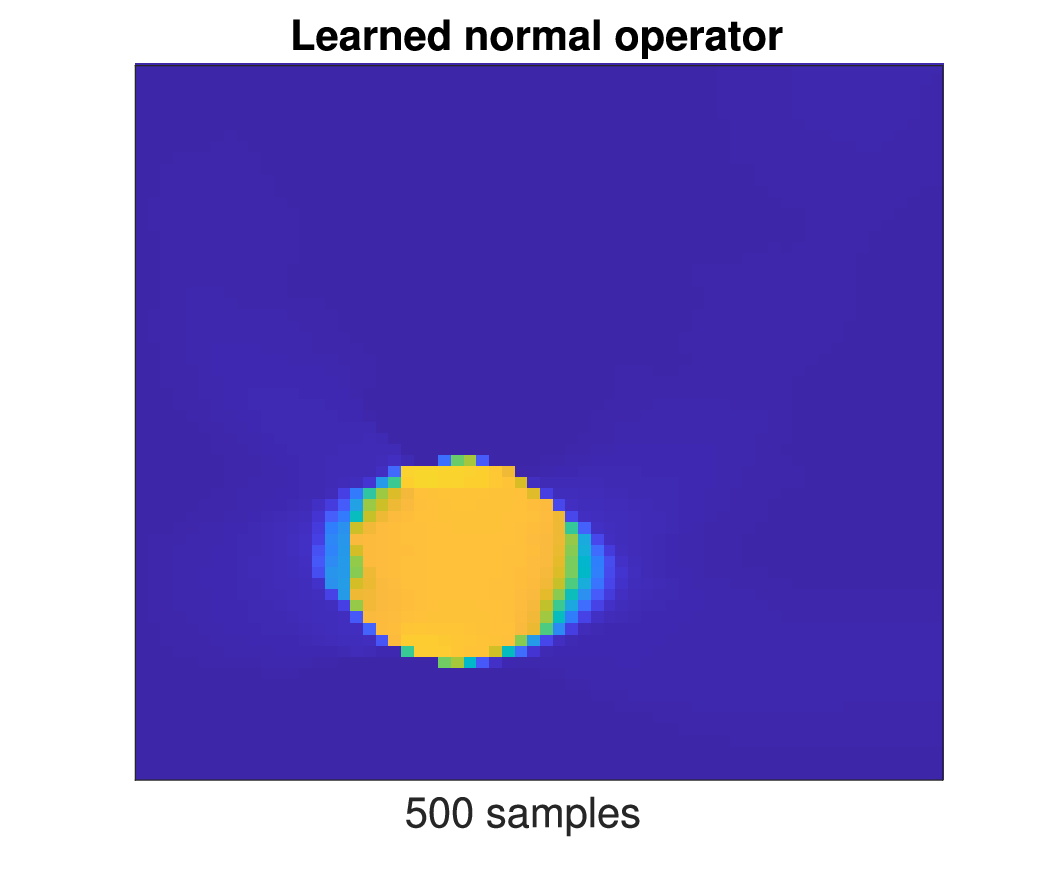}};
    \node[above=0 of image] {\scriptsize Learned normal operator};
    \end{tikzpicture}
    \end{subfigure} \\[5pt]
    ~~~~~~\begin{subfigure}[t]{0.4\textwidth}
    \begin{tikzpicture}
    \node[inner sep=0] (image) at (0,0) {\adjincludegraphics[width=\textwidth,trim={0 0 0 {.07\height}},clip]{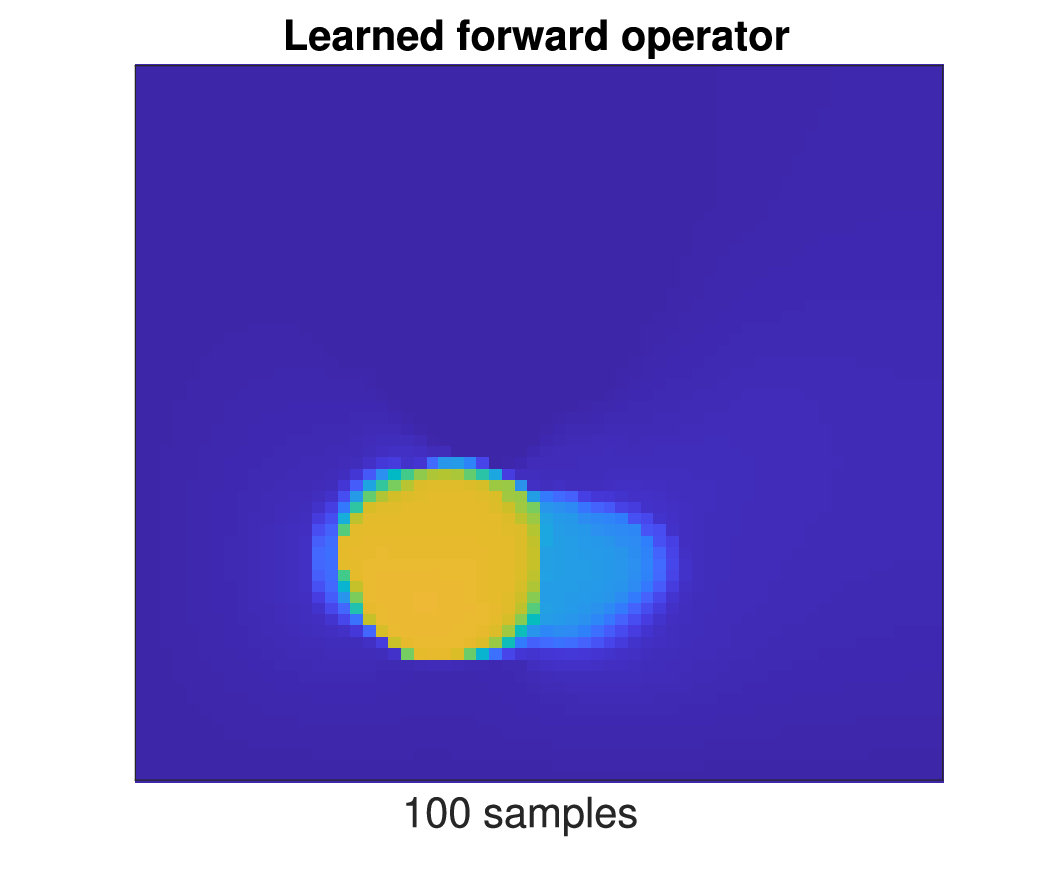}};
    \node[above=0 of image] {\scriptsize Learned forward operator};
    \end{tikzpicture}
	\end{subfigure}
	\begin{subfigure}[t]{0.4\textwidth}
    \begin{tikzpicture}
    \node[inner sep=0] (image) at (0,0) {\adjincludegraphics[width=\textwidth,trim={0 0 0 {.07\height}},clip]{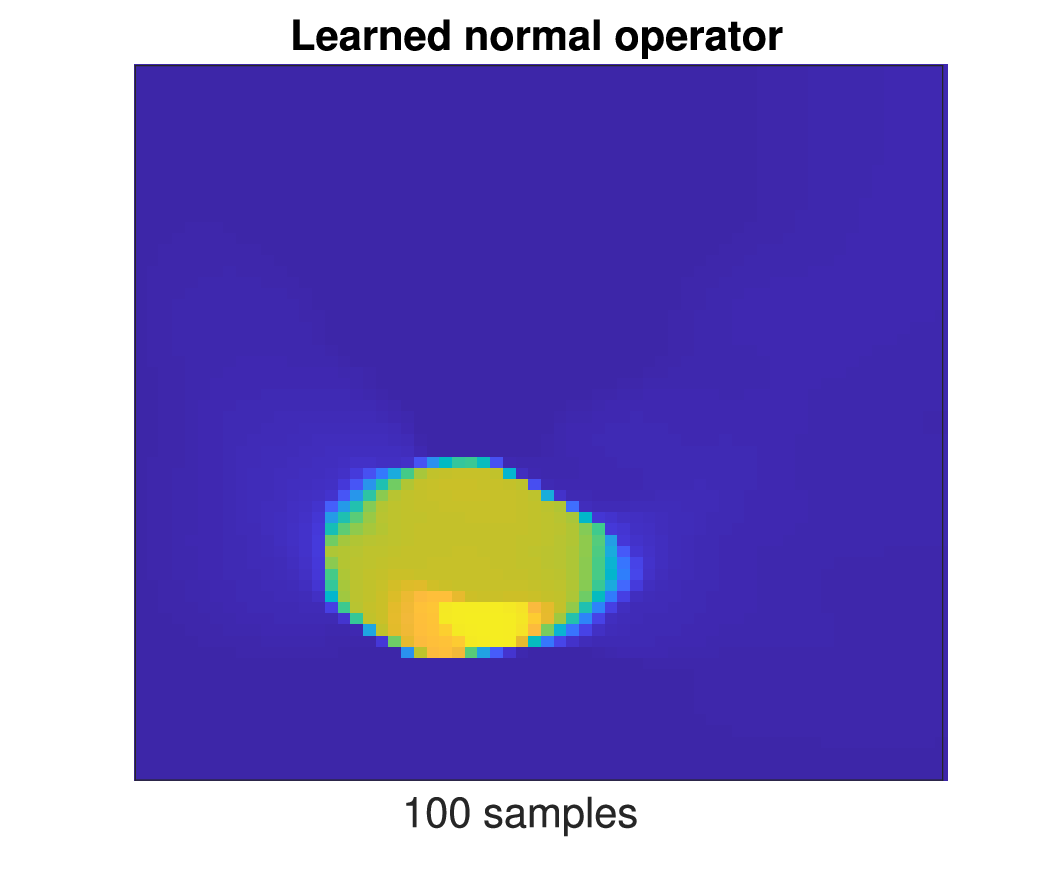}};
    \node[above=0 of image] {\scriptsize Learned normal operator};
    \end{tikzpicture}
    \end{subfigure} \\[5pt]
    ~~~~~~\begin{subfigure}[t]{0.4\textwidth}
    \begin{tikzpicture}
    \node[inner sep=0] (image) at (0,0) {\adjincludegraphics[width=\textwidth,trim={0 0 0 {.07\height}},clip]{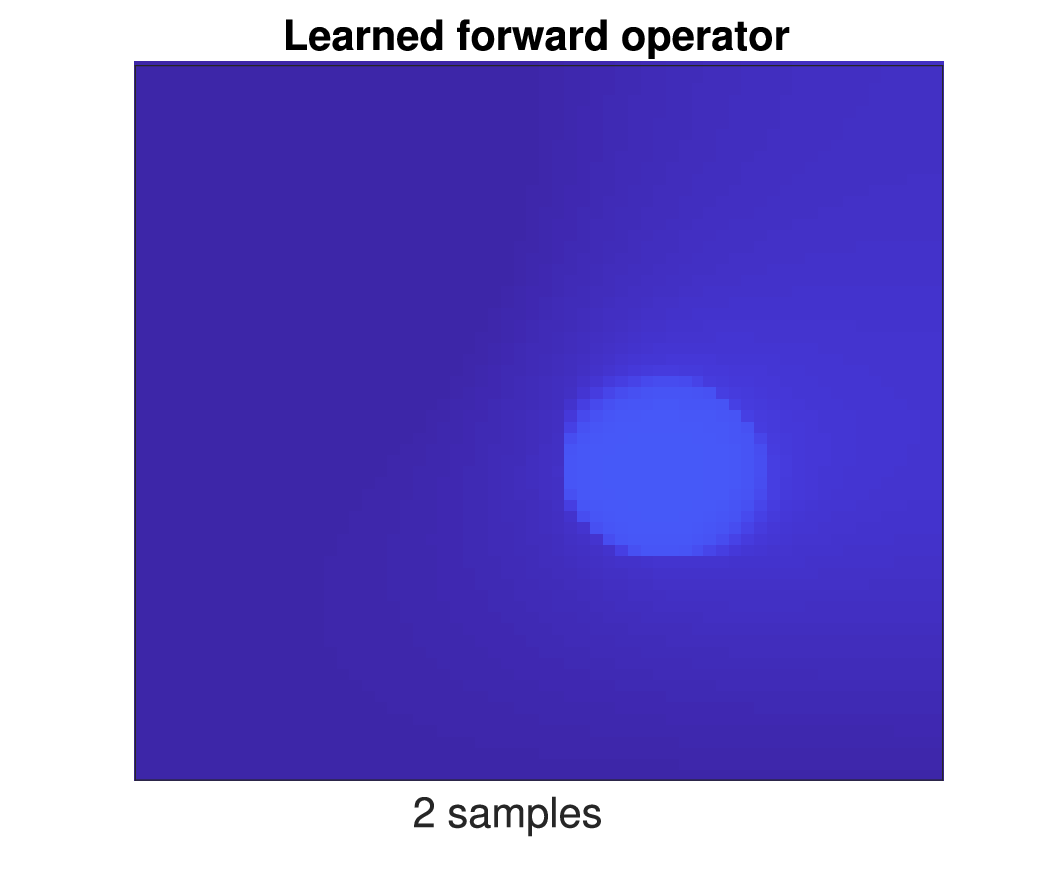}};
    \node[above=0 of image] {\scriptsize Learned forward operator};
    \end{tikzpicture}
	\end{subfigure}
	\begin{subfigure}[t]{0.4\textwidth}
    \begin{tikzpicture}
    \node[inner sep=0] (image) at (0,0) {\adjincludegraphics[width=\textwidth,trim={0 0 0 {.07\height}},clip]{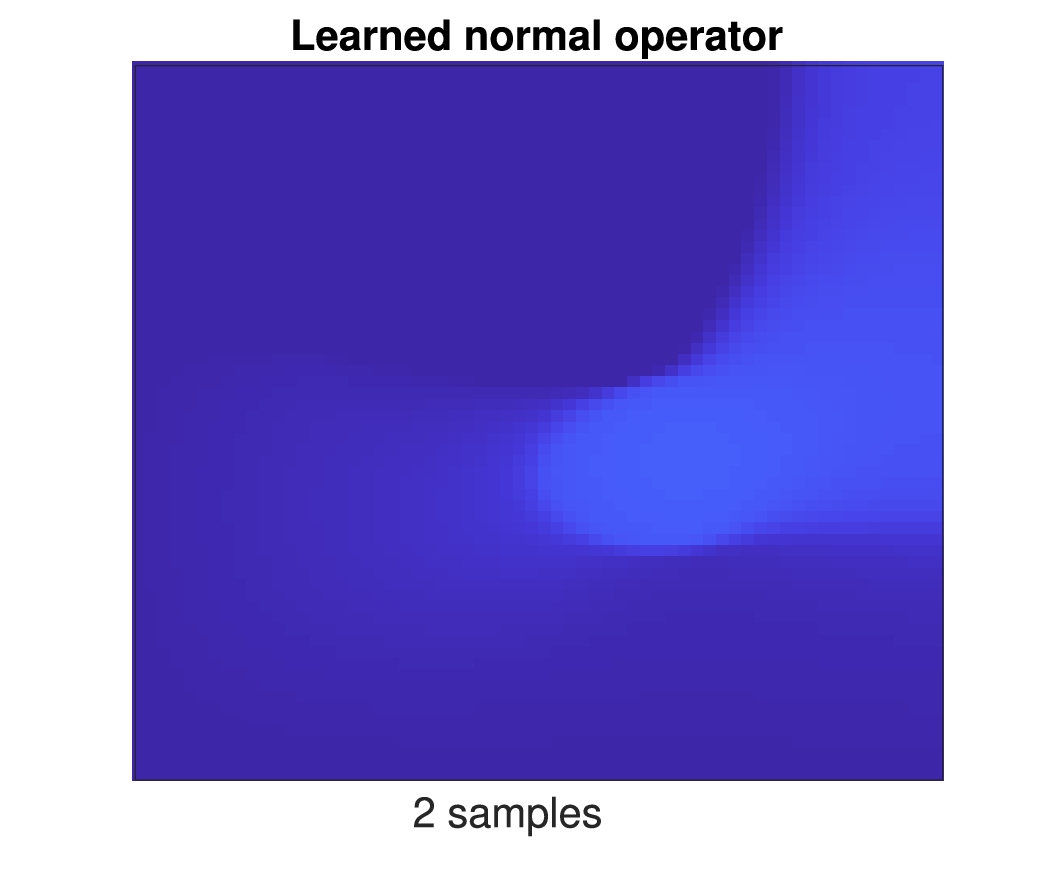}};
    \node[above=0 of image] {\scriptsize Learned normal operator};
    \end{tikzpicture}
    \end{subfigure}
    \end{minipage}%
    \begin{minipage}{0.2\textwidth}
        \includegraphics[width=0.5\textwidth]{fig_proj/colour-bar.png}
    \end{minipage}
    \caption{Projected variational regularisation. The effect of the number of training samples. The method is surprisingly robust with respect to the number of samples. Reconstructions are still reasonable for as few as $500$ samples (top row), and even for $100$ (middle row) samples the location of the disc is identified reasonably well.  If we decrease the number of samples further so that the support of the samples does not overlap with the support of the ground truth then the location of the disc is identified incorrectly (bottom row). }\label{fig:proj-var-number-of-samples}
\end{figure}

Using a combined model as in \cref{eq:var-reg-GD-combined,eq:var-reg-GD-normal-combined}, surprisingly, yields inferior performance. For the combined forward operator (bottom row left) we see artefacts outside of the support of the discs in the training set; this is the region where only the approximate model is available. The learned normal operator (bottom row right) is affected less, possibly because the support of the backprojections $\trback$ is larger than that of the training images $\trim$. However, the reconstruction still looks blurrier than without the approximate model (middle row right).

In \cref{fig:proj-var-number-of-samples} we investigate the influence of the size of the training set $n$ on the reconstruction quality. We find that the method is very robust. For as few as $n=500$ samples we still obtain a reasonable reconstruction (top row), the learned normal operator (right) performing slightly better than the learned forward operator (left). Even for $n=100$ the location of the disc is identified reasonably well, even if the shape is not well reconstructed  (middle row; the regularisation parameter has to be chosen larger in this case to ensure stability). The reason for such robustness seems to be that the learned model ``sees'' the area where the ground truth disc is located because there are other discs nearby in the training set. If we further deflate the training set  to an extremely low size of $n=2$ so that there is no overlap between the training discs and the exact solution, then the support of the reconstruction is completely off (bottom row).

%%%%%%%%%%%%%%%%%%%%%%%%%%%%%%%%%%%%%%%%%%%%%%%%
\subsection{Experiments with learned model correction}
We continue with experiments with the second approach considered in this chapter, the learned model correction. 
As we have seen in Figure \ref{fig:ForwApprox}, the approximate model  introduces artefacts that will cause the gradients to be incorrect,  and hence a correction needs to be applied.
As discussed in \cref{sec:model-correction}, there are several ways to achieve such a model correction, classified into implicit and explicit approaches. Here, we will only discuss the explicit corrections.

\begin{figure}[t!]
    \centering
        \includegraphics[width=.8\textwidth]{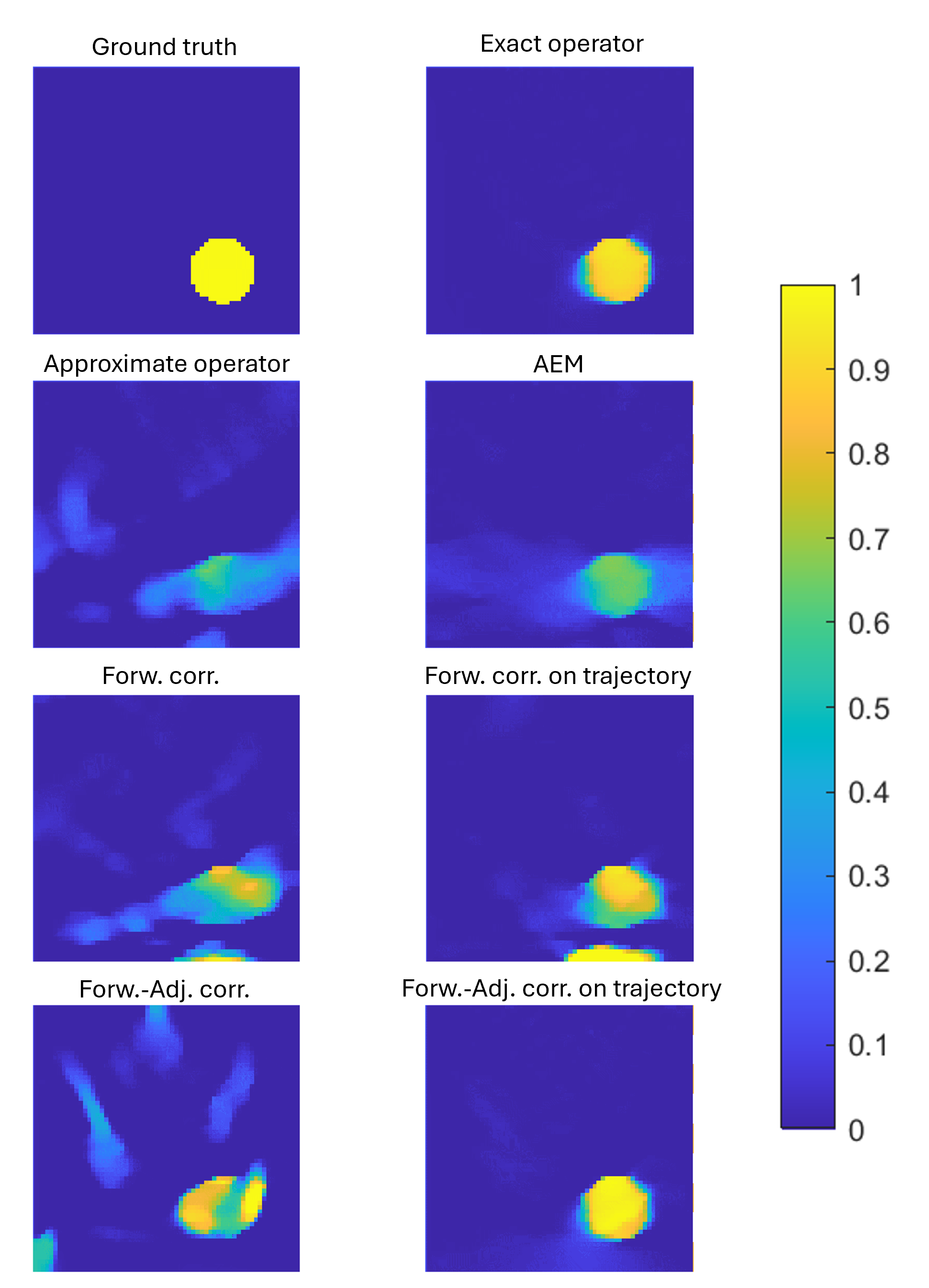}
        \caption{\label{fig:model_corr} Comparison of reconstructions using the approximate model and learned corrections under various training schemes. (Top row) Ground truth and reconstruction with the accurate model, the sensor is located at the top. (Second row) Reconstruction using the approximate operator $\widetilde{A}$ and reconstruction using the approximation error method. (Third row) Learned correction of the forward model trained on disks only (left) and recursively on the trajectory (right). (Bottom row) Forward-adjoint correction trained on disks only (left) and recursively on the trajectory (right). }
\end{figure}

\begin{figure}[t!]
    \centering
        \includegraphics[height=.45\textwidth]{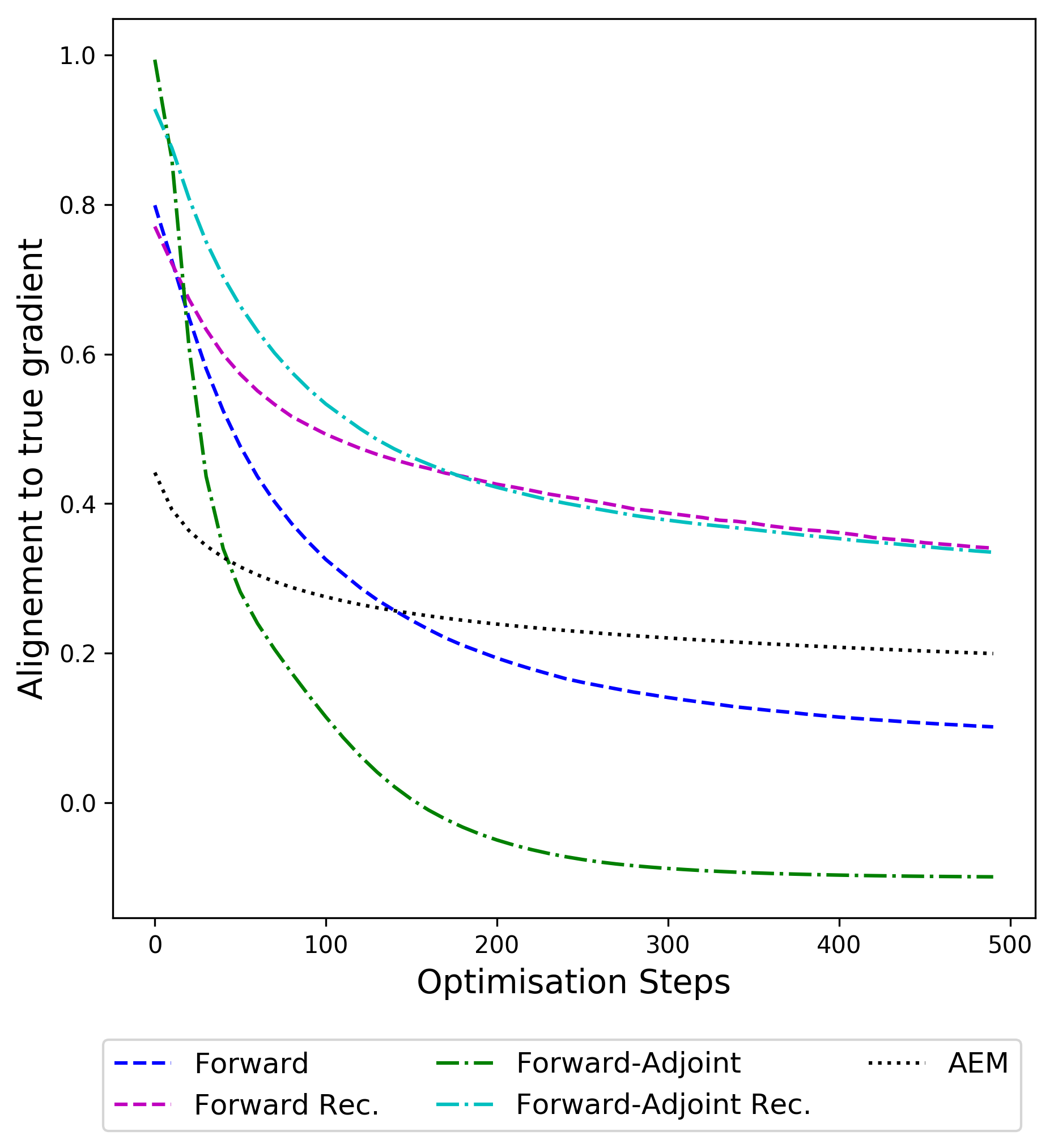}
        \,
        \includegraphics[height=.45\textwidth]{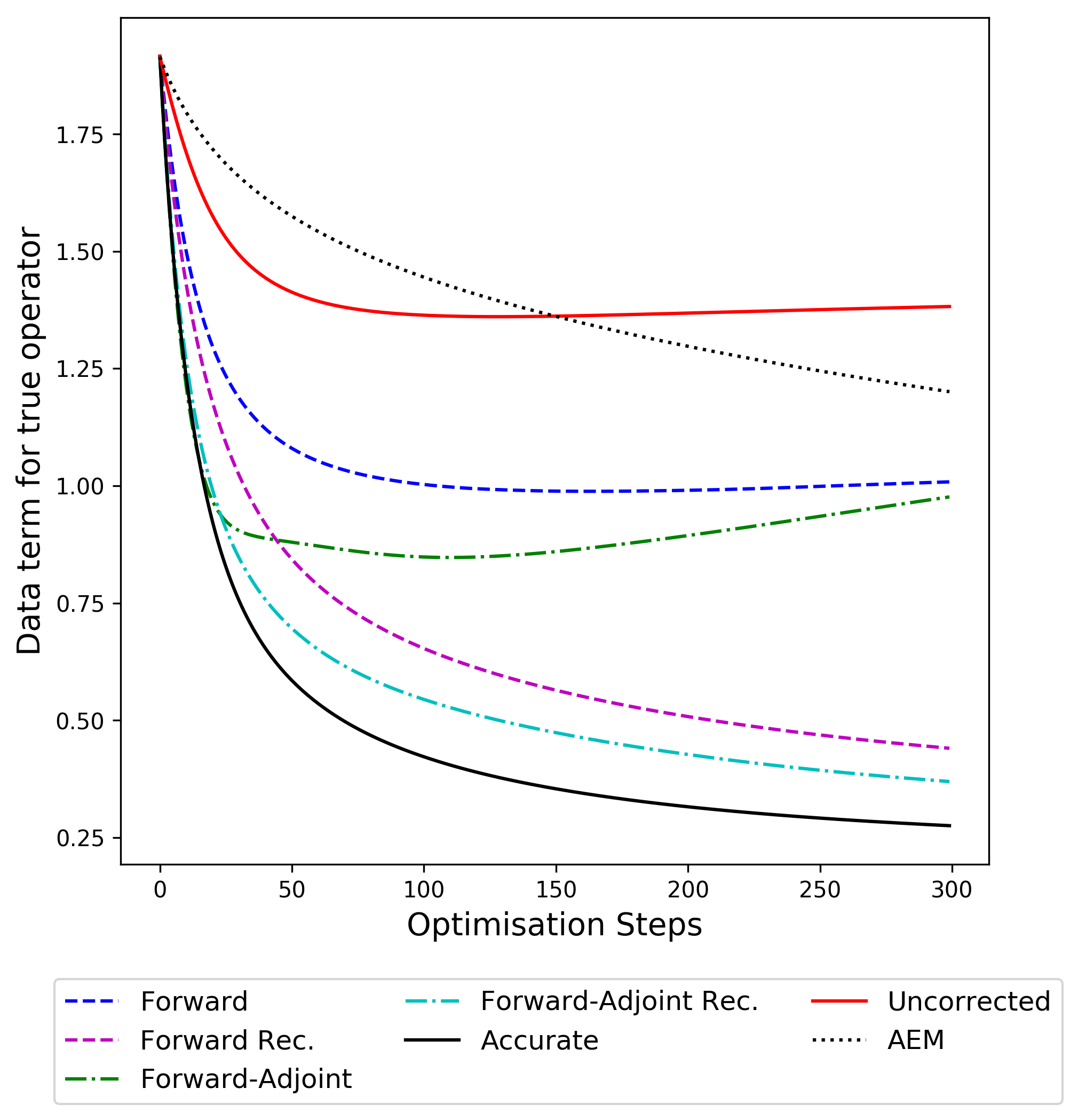}
        \caption{(Left) Alignment \eqref{eqn:AligenementVariational} of the approximate gradient and the gradient of the accurate data term $\TrueOp^*(\TrueOp \iter{x}{k} - y)$ for different training approaches. (Right)
        True data term $\|A\iter{x}{k} - y\|_{\calY}$  for different methods, tracked throughout the gradient descent scheme.}
        \label{fig:Conv_plots}
\end{figure}

The experiments are also conducted on the disc data set and we compare learning a correction for the forward operator only with corrections for both the forward operator and its  adjoint. Additionally, we discuss the influence of training only on the data manifold of ground-truth samples $\trim$ or along the trajectory (recursive training) as discussed in \cref{sec:forw-adj-cor}. Experiments presented here were first reported in \cite{lunz21}. The reconstructions are shown in \cref{fig:model_corr} and the convergence plots in \cref{fig:Conv_plots}.

In \cref{fig:model_corr} we see that the accurate model with total variation regularisation---our reference solution---produces a reasonably good reconstruction with an average error of 12\% over the 64 test samples, but due to strong limited view artefacts we still have a slight smearing visible. 
The approximate and uncorrected forward operator $\widetilde{A}$ is not able to produce a sufficiently good reconstruction and  causes  strong artefacts in the background medium. This is also reflected in the average reconstruction error of 55\%. The classical approximation error method is indeed able to correct the strong artefacts and reaches a relative error of 32\%, but results in a loss of contrast and thus wrong quantitative values. It can be also seen in Figure \ref{fig:Conv_plots} that the convergence is much slower.

Let us now examine the learned model corrections using a neural network, that is  $\CorrectedOp := \ForwardCor \circ \ApproxOp$ for the forward correction and, if used, $\CorrectedAd := \AdjointCor \circ \ApproxOp^*$ for the adjoint correction. We can see a clear difference between training along the trajectory compared to training on the data manifold only. For both corrections, forward and forward-adjoint, the training on the simple data manifold is not sufficient and leads to strong artefacts, resulting in an average relative error of 53\% and 41\% respectively. This is due to the correction not being valid near the point of convergence and hence the conditions of Theorem \ref{thm:conv-neighbourhood} are violated (recall that the theorem requires that the correction  is valid for all iterates $\iter{x}{k}$). If we make sure that the correction is valid for all iterates by training the networks with \eqref{equ:forwardLoss_FBC} not only on the data manifold of ground truth samples, but also for all $\iter{x}{k}$ that arise during the optimisation procedure, then we can ensure convergence to a neighbourhood of the accurate minimiser as stated by Theorem \ref{thm:conv-neighbourhood} for the forward-adjoint correction. This is clearly seen in the reconstruction that is visually close to the accurate model and also reflected in the average relative error of 14\% closest to the accurate model.

The effect of gradient alignment can also be observed in the convergence plots in \cref{fig:Conv_plots}. Recursive training, both for the forward and forward-adjoint corrections, ensures positive alignment between the true and approximate gradients and results in convergence of the data misfit (measured using the accurate model). Without recursive training, both forward and forward-adjoint corrections result in non-convergence of the data misfit.

Unfortunately, this improvement in reconstruction quality requires  an expensive training procedure that takes about 4 days in total to train the networks carefully, whereas training on the data manifold of discs and corresponding measurements takes only a few hours for the forward-adjoint correction. This indicates that there is a need to improve training strategies as we discuss next.

%%%%%%%

\subsubsection{Training without trajectory} \label{sec:training-trajectory}
We have seen in the previous section that  training along a trajectory is crucial for the success of the explicit learned model correction when solving the variational problem, which is in accordance with Theorem \ref{thm:conv-neighbourhood} requiring a valid approximation for all iterates $\iter{x}{k}$ obtained during minimisation.
The high computational cost of such recursive training makes it necessary to think about alternatives. One possibility is to restrict the minimisation of the variational problem 
\begin{equation}\label{equ:variationalSolutionEquality}
    \argmin_{x\in \calX} \frac{1}{2} \|\CorrectedOp(x) - y\|_{\calY}^2 + \lambda \RegularisationOp(x)
\end{equation}
to a suitable set, such as the manifold representing the data distribution $\trim$. This way, the correction can be trained just on the manifold itself and constraining the trajectory to the vicinity of the manifold will ensure that the correction is valid for all $\iter{x}{k}$. 

\begin{figure}[!ht]
    \centering
        \includegraphics[width=.6\textwidth]{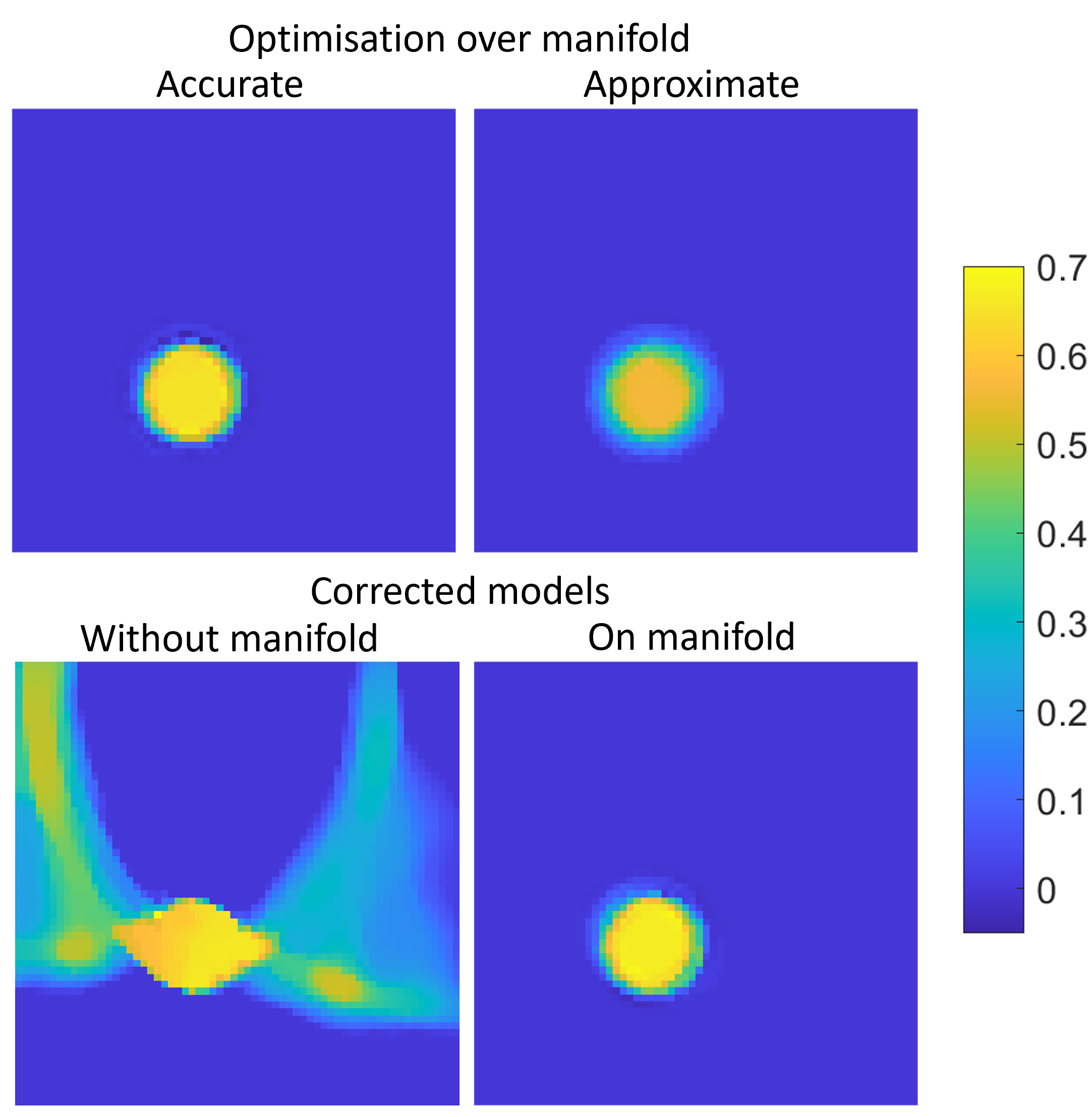}
        \caption{\label{fig:manifold} Proof-of-concept results for training and optimisation over the data manifold instead of the trajectory. (Top) optimisation is carried out on the data manifold. This already improves reconstructions with the approximate model, but results in loss of contrast. (Bottom) Corrected models trained on data manifold only. On the left optimisation is carried out over all $\calX$ and on the right only on the data manifold. }
\end{figure}

We present, in Figure \ref{fig:manifold}, a proof-of-concept result for the hypothesis that training and optimisation over the data manifold can eliminate the need for trajectory training. We see that solving the variational problem on the manifold  works well for the accurate model, whereas the approximate model suffers loss of contrast and sharpness. If we use the corrected model trained on the data manifold, but optimisation is performed freely in the full space, we  introduce strong artefacts, as the correction is not valid for all iterates. However, when the optimisation path is restricted to the manifold we obtain a result close to that of the accurate model. Training of the model correction on the data manifold only takes roughly 90 minutes, compared to the full trajectory training that requires 4 days. This is a promising solution to the trajectory training problem and is currently work in progress.

%%%%%%%%%%%%%%%%%%%%%%%%%%%%%%%%%%%%%%%%%%%%%%%%
\section{Summary and conclusions}
\label{sec:conclusions}
Inverse problems involving high dimensions and/or computationally expensive forward operators necessitate the use of computationally cheaper approximate models. This applies, e.g., to settings when imaging is performed in time-critical scenarios or when the forward operator is called many times within a larger pipeline of a learning framework. 

In this chapter we   have discussed two data-driven approaches, learning an approximation of the forward or the inverse operator using data-driven projections \cite{asp-kor-sch-2020}, and a data-driven correction to an analytic approximation \cite{lunz21}. A common theme has emerged, which is that in addition to learning an approximation of the forward operator, one often needs to learn a separate approximation of the adjoint.

While both approaches provide a possible solution to the problem of computationally expensive models in inverse problems, they also come with some drawbacks. The data-driven projection method requires a good quality \textit{linear} approximation of the ground truth by training images, hence, for example, it is sensitive to shifting the image. If the forward operator is shift-equivariant (e.g., a convolution), it could be possible to incorporate a non-linear ``projection'' onto the training set by finding the shift that minimises the distance between the image and the span of the training data (a similar approach can be applied, e.g., to rotations). However,  this is a research direction not yet taken.

In the case of learned model corrections we have discussed that while solutions are faster to compute, the computational burden moves to an expensive training phase that needs to ensure validity of the correction for all iterates of the optimisation trajectory to ensure convergence. Finally, we presented a proof-of-concept solution that may overcome this problem by limiting training and optimisation to the data manifold, which is an interesting direction for future studies.

%%%%%%%%%%%%%%%%%%%%%%%%%%%%%%%%%%%%%%%%%%%%%%%%
\section*{Acknowledgements}
%%%%%%%%%%%%%%%%%%%%%%%%%

SA and YK acknowledge EPSRC grant EP/V026259/1. YK also acknowledges the EPSRC Fellowship EP/V003615/2 and the support of the National Physical Laboratory. AH acknowledges support by the Research Council of Finland: Academy Research Fellowship (Project 338408), the Centre of Excellence of Inverse Modelling and Imaging project (Project 353093), and the Flagship of Advanced Mathematics for Sensing Imaging and Modelling grant (Project 359186).

The authors would also like to thank the Isaac Newton Institute for Mathematical Sciences, Cambridge, for support and hospitality during the programme Mathematics of deep learning (supported by EPSRC grant EP/R014604/1).

\printbibliography

\end{document}